\documentclass[11pt]{article}
\usepackage{amssymb,amsmath,latexsym}
\usepackage{epsfig}

\usepackage{epic,eepic, color}

\newtheorem{thm}{Theorem}[section]

\newtheorem{lemma}[thm]{Lemma}
\newtheorem{prop}[thm]{Proposition}

\numberwithin{equation}{section}

\def\pf{{\medskip\noindent {\bf Proof. }}}
\def\qed{{\hfill $\Box$ \bigskip}}
\def\R{{\mathbb R}}
\def\P{{\mathbb P}}
\def\E{{\mathbb E}}
\def\1{{\bf 1}}

\renewcommand{\Pr}{\mathbb P}

\def\bP {{\mathbb P}}  \def\bR {{\mathbb R}}

\def\R {{\mathbb R}}

\def\nn{\nonumber}

\def\wt{\widetilde}
\def\wh{\widehat}

\def\E{{\mathbb E}}
\def\P{{\mathbb P}}

\def\bea{\begin{align*}}
\def\eea{\end{align*}}
\def\bee{\begin{equation}}
\def\eee{\end{equation}}

\def\eps{\varepsilon}

\def\wh{\widehat}

\makeatletter
\@addtoreset{equation}{section}

\makeatother

\begin{document}
\allowdisplaybreaks
\bibliographystyle{plain}

\title{\Large \bf
Global Heat Kernel Estimate for Relativistic Stable Processes \\
in Exterior Open Sets}

\author{{\bf Zhen-Qing Chen}\thanks{Research partially supported
by NSF Grants DMS-0906743 and DMR-1035196.}, \quad {\bf Panki Kim}\thanks{This research was supported by Basic Science Research Program through the National Research Foundation of Korea (NRF) funded by the Ministry of Education, Science and Technology(0409-20110087). } \quad and  \quad {\bf Renming Song}\thanks{Research supported in part by a grant from the Simons Foundation (208236).}
}
\date{(December 9, 2011)}

\maketitle

\begin{abstract}
In this paper, sharp
two-sided estimates for the transition densities of
relativistic $\alpha$-stable processes with mass $m\in (0, 1]$ in $C^{1,1}$ exterior open sets
are established for all time $t>0$.
These transition densities are also the Dirichlet heat kernels of
$m-(m^{2/\alpha}-\Delta)^{\alpha/2}$ with $m\in (0, 1]$ in $C^{1,1}$ exterior open sets.
The estimates are uniform in $m$ in the sense that the
constants are independent of $m\in (0, 1]$.
As a corollary of our main result, we establish
sharp two-sided Green function estimates for relativistic $\alpha$-stable processes with mass $m\in (0, 1]$ in $C^{1,1}$ exterior open sets.
\end{abstract}


\bigskip
\section{Introduction}

Let $d\ge 1$ and $\alpha\in (0, 2)$. For any $m\ge 0$, a relativistic
$\alpha$-stable process $X^m$ in $\bR^d$ with mass $m$ is a L\'evy process
with characteristic function given by
 \bee\label{e:ch}
   \E \left[\exp\left(i \xi \cdot (X^m_t-X^m_0)
   \right) \right] = \exp\left(-t \left( \big(|\xi|^2+ m^{2/\alpha}
    \big)^{\alpha/2}-m\right) \right), \qquad \xi \in \bR^d.
 \eee
When $m=0$, $X^m$ is simply a (rotationally) symmetric
$\alpha$-stable process in $\bR^d$. The infinitesimal generator of
$X^m$ is $m-(-\Delta +m^{2/\alpha} )^{\alpha /2}$.
 When $\alpha=1$, the infinitesimal generator reduces
to the free relativistic Hamiltonian $ m - \sqrt{-\Delta + m^{2}}$.
There exists a huge literature on the properties of
relativistic Hamiltonians (for example, see
\cite{CMS,FL,He,Lieb,LY}). Relativistic $\alpha$-stable processes
have been studied recently in \cite{CK2, CS4, GR2, K2, KL, KS, R}.

Recall that an open set $D$ in $\bR^d$ (when $d\ge 2$) is said to be
a (global) $C^{1,1}$ open set if there exist a localization radius $ r_0>0 $
and a constant $\Lambda_0>0$ such that for every $z\in\partial D$,
there exist a $C^{1,1}$-function $\phi=\phi_z: \bR^{d-1}\to \bR$
satisfying $\phi (0)=0$, $\nabla\phi (0)=(0, \dots, 0)$, $\| \nabla
\phi  \|_\infty \leq \Lambda_0$, $| \nabla \phi (x)-\nabla \phi (z)|
\leq \Lambda_0 |x-z|$, and an orthonormal coordinate system $y=(y_1,
\cdots, y_{d-1}, y_d):=(\wt y, \, y_d)$ such  that $ B(z, r_0 )\cap
D=B(z, r_0 )\cap \{ y: y_d > \phi (\wt y) \}$.
We call the pair $(r_0, \Lambda_0)$ the
characteristics of the $C^{1,1}$ open set $D$.
By a $C^{1,1}$ open set
in $\bR$ we mean an open set which can be expressed
 as the union of disjoint intervals so that the minimum of
the lengths of all these intervals is positive and the minimum of
the distances between these intervals is positive. Note that a
$C^{1,1}$ open set can be unbounded and disconnected.

For an open set $D\subset \R^d$, let $X^{m, D}$ be the subprocess
of $X^m$ killed upon exiting $D$. It is easy to see
(cf. \cite{CK2}) that $X^{m, D}$ has a jointly continuous
transition density function
$p^m_D(t, x, y)$ with respect to
the Lebesgue measure on $D$.
$p^m_D$ is also called the Dirichlet heat kernel of $m-(-\Delta +m^{2/\alpha} )^{\alpha /2}|_D$
with zero exterior condition.

A relativistic $\alpha$-stable process is a discontinuous Markov process.
Sharp estimates on the transition density functions of discontinuous
Markov processes are of current research interests
(see \cite{BBCK, C, CKK, CKK2, CKK3, CK, CK2} and the reference
therein). Dirichlet heat kernel estimates for symmetric  stable
processes were first obtained in \cite{CKS} on $C^{1,1}$ open sets
for $t\leq 1$ and for all $t>0$ when the $C^{1,1}$ open set is bounded.
In \cite{BGR}, Dirichlet heat kernel estimates for symmetric  stable
processes for a large class of non-smooth open sets were obtained in terms of surviving probabilities.
In \cite{CT}, global Dirichlet heat kernel estimates for symmetric
stable processes are derived for $C^{1,1}$  exterior open sets
as well as  for half-space-like open sets.
The ideas of \cite{CKS} have been adapted to establish  sharp two-sided estimates
for the Dirichlet heat kernels of other discontinuous Markov processes in open sets,
see \cite{CKS1, CKS3, CKS5}.
In particular, the following result is established
 in \cite[Theorem 1.1]{CKS3}.
In this paper, for any $a, b\in \bR$, we use the notations
$a\wedge b:=\min \{a, b\}$ and $a\vee b:=\max\{a, b\}$.

\begin{thm}\label{t:main}
Suppose that
$\alpha\in (0, 2)$ and
$D$ is a $C^{1,1}$ open set in $\bR^d$ with
$C^{1,1}$ characteristics $(r_0, \Lambda_0)$.
Let $\delta_D(x)$
be the Euclidean distance between $x$ and $D^c$.
\begin{description}
\item{\rm (i)}
For any $M>0$ and $T>0$,
there exists $c_1=
c_1(d, \alpha, r_0, \Lambda_0, M, T)>1$
such that for any $m\in (0, M]$ and
$(t,x,y) \in (0, T]\times D\times D$,
  \begin{align}
& \hskip 0.3truein
   \frac1{c_1} \Big( 1\wedge
\frac{\delta_D(x)^{\alpha/2}}{\sqrt{t}}\Big) \Big( 1\wedge
\frac{\delta_D(y)^{\alpha/2}}{\sqrt{t}}\Big) \Big( t^{-d/\alpha}
\wedge \frac{t
  \phi(m^{1/\alpha} |x-y|)}
{|x-y|^{d+\alpha}}\Big) \nonumber\\
&\hskip 0.3truein \le p^m_D(t, x, y)\nonumber\\
&\hskip 0.3truein  \le c_1  \Big( 1\wedge
\frac{\delta_D(x)^{\alpha/2}}{\sqrt{t}}\Big) \Big( 1\wedge
\frac{\delta_D(y)^{\alpha/2}}{\sqrt{t}}\Big) \Big( t^{-d/\alpha}
\wedge \frac{t
\phi(m^{1/\alpha} |x-y|/
(16))}
{|x-y|^{d+\alpha}}\Big)
, \label{e:1.2}
 \end{align}
where $\phi (r)= e^{-r}(1+r^{(d+\alpha-1)/2})$.

\item{\rm (ii)} Suppose in addition that $D$ is bounded.
 For any $M>0$ and $T>0$,
there exists $c_2=
c_2(d, \alpha, r_0, \Lambda_0, M, T, \text{diam}(D)) >1$
such that for any $m\in (0, M]$ and
$(t,x,y) \in [T, \infty)\times D\times D$,
\begin{eqnarray*}
c_2^{-1}\, e^{- t\, \lambda^{\alpha, m, D}_1  }\, \delta_D
(x)^{\alpha/2}\, \delta_D (y)^{\alpha/2}
\leq
p^m_D(t, x, y)
 \leq
c_2\, e^{-t \, \lambda^{\alpha, m, D}_1  }\, \delta_D
(x)^{\alpha/2} \,\delta_D (y)^{\alpha/2} ,
\end{eqnarray*}
where $\lambda^{\alpha, m, D}_1>0$ is the smallest eigenvalue of the
restriction of $(m^{2/\alpha}-\Delta)^{\alpha /2}-m$ in $D$ with
zero exterior condition.
\end{description}
\end{thm}

Note that, although the small time estimates on $p^m_D(t, x, y)$ in Theorem
\ref{t:main}(i) are valid for all $C^{1,1}$ open sets, the large time estimates
in Theorem \ref{t:main}(ii) are only for bounded $C^{1,1}$ open sets.
As one sees
for the case of symmetric $\alpha$-stable processes in \cite{CT},
the large time heat kernel estimates for unbounded open sets
are typically very different from that in the bounded open sets and depend
on the geometry of the unbounded open sets.
Sharp two-sided estimates on $p^m_D(t, x, y)$ valid for all time $t>0$
have recently been established for half-space-like $C^{1,1}$ open sets
in \cite{CKS5} by using some ideas from \cite{CT}. The goal of this paper is
to establish sharp  two-sided estimates on $p^m_D(t, x, y)$ for exterior
$C^{1,1}$ open sets
that hold
for all $t>0$.

Recall that an open set $D$ in $\R^d$ is called an exterior open set if $D^c$ is compact.
For any $m, b, c>0$, we define a function
$\Psi_{d, \alpha, m, b, c} (t, x, y)$
 on $(0, \infty)\times \R^d\times \R^d$ by
\begin{align}
\Psi_{d, \alpha, m, b, c} (t, x, y)
:= \begin{cases}
 {t^{-d/\alpha}}\wedge\frac{t
  \phi(c^{-1}m^{1/\alpha} |x-y|)} {|x-y|^{d+\alpha}}
  ~ &\hbox{when }   t \in (0,b/m],
\\
 m^{d/\alpha-d/2}t^{-d /2}\exp\left( -c^{-1} (m^{1/\alpha}
|x-y|\wedge m^{2/\alpha-1}\frac{|x-y|^2}{t}   )\right)
~ &\hbox{when }
t \in (b/m, \infty),
\end{cases}\label{e:Psi}
\end{align}
where $\phi (r)= e^{-r} \left(1+r^{(d+\alpha-1)/2}\right)$.
The following is the main result of this paper.

\begin{thm}
\label{t:MAIN}
Suppose that
$\alpha\in (0, 2)$,
$d\ge 3$, $M>0$, $b>0$, $R>0$ and $D$ is an exterior $C^{1,1}$ open set in
$\R^d$ with $C^{1,1}$
characteristics $(r_0, \Lambda_0)$ and $D^c\subset B(0, R)$.
Then there are constants
$c_i=c_i(d, \alpha,  M, b, r_0, \Lambda_0, R)>1, i=1,2,$ such that
for every $m\in (0, M
]$,  $t>0$ and
$(x, y)\in D \times D$,
$$
p^m_D(t,x,y)\le c_1\left( 1\wedge \frac{\delta_D(x)}{ 1 \wedge t^{1/\alpha}} \right)^{\alpha/2}
\left( 1\wedge \frac{\delta_D(y)}{1 \wedge  t^{1/\alpha}}\right)^{\alpha/2}
\Psi_{d, \alpha, m, b, c_2} ( t, x, y)
$$
and
$$
p^m_D(t,x,y)\ge c_1^{-1} \left( 1\wedge \frac{\delta_D(x)}{ 1 \wedge t^{1/\alpha}} \right)^{\alpha/2}
\left( 1\wedge \frac{\delta_D(y)}{1 \wedge  t^{1/\alpha}}
\right)^{\alpha/2} \Psi_{d, \alpha, m, b,  1/c_2} ( t, x, y).
$$
\end{thm}

\medskip
It is known (see Theorem \ref{T1u} below) that
there are constants $c_3>1$ and $c_4\geq 1$ such that for all $m>0$
and $(t, x, y)\in (0, \infty) \times \bR^d\times \bR^d$,
\begin{equation}\label{e:newer1}
  c_3^{-1}  \Psi_{d, \alpha, m, 1,  1/c_4} ( t, x, y)
\leq p^m(t, x, y) \leq c_3 \Psi_{d, \alpha, m, 1,  c_4} ( t, x, y).
\end{equation}
By integrating the sharp heat kernel estimates in Theorem \ref{t:MAIN} (with $b=1$)
over $y\in D$ and using \eqref{e:newer1},
 one can easily conclude that there is a constant
$c_5=c_5 (d, \alpha,  M, r_0, \Lambda_0, R)\geq 1$ so that
for every $m \in (0, M]$, $x\in D$ and $t>0$,
$$ c_5^{-1} \left( 1\wedge \frac{\delta_D(x)}{ 1 \wedge t^{1/\alpha}} \right)^{\alpha/2} \leq \P_x (\tau^m_D >t) \leq c_5
 \left( 1\wedge \frac{\delta_D(x)}{ 1 \wedge t^{1/\alpha}} \right)^{\alpha/2},
 $$
 where $\tau^m_D=\inf\{ t>0: X^m_t \notin D\}$.
We emphasize that the sharp heat kernel estimates in Theorem \ref{t:MAIN}
hold uniformly in $m\in (0, M]$. Thus passing $m\downarrow 0$
recovers the sharp heat kernel
estimates for symmetric $\alpha$-stable processes in exterior
$C^{1,1}$ open sets
when dimension $d\geq 3$ that were previously obtained in \cite{CT}.
(The estimates in \cite{CT} hold for every $d\geq 2$.)
The large time upper bound estimate in Theorem \ref{t:MAIN}
is quite easy to establish, which is given
at the end of Section \ref{S:2}.
The main task of this paper is to establish the large
time lower bound estimate for $p^m_D(t, x, y)$.
Comparing with the case of symmetric stable processes,
due to the fact that the L\'evy densities of
relativistic stable processes decay exponentially fast at infinity,
the large time lower bound estimates
 for $p^m_D$ is much harder to establish.
The reason that we assume $d\geq 3$ in Theorem \ref{t:MAIN} is
that, due to Chung-Fuck's recurrence criterion for L\'evy processes,
relativistic stable processes are transient if and only if
$d\geq 3$.

Integrating the heat kernel estimates in Theorem \ref{t:MAIN}
in $t\in (0, \infty)$,
one gets the following sharp two-sided Green function estimates of $X^m$
in exterior $C^{1,1}$ open
sets, which is uniform in $m\in (0, M]$.

\begin{thm}\label{T:1.3}
Suppose that $d\ge 3$, $M>0, R>0$ and $D$ is an exterior $C^{1,1}$ open set in
$\R^d$ with $C^{1,1}$ characteristics $(r_0, \Lambda_0)$ and   $D^c\subset B(0, R)$.
Then there is a constant
$c=c(d, \alpha,  M, r_0, \Lambda_0, R)>1$ such that
for every $m\in (0, M
]$ and
$(x, y)\in D \times D$,
\begin{align*}
&c^{-1} \frac{1 +
 (m^{1/\alpha} |x-y|)^{2-\alpha}     }{|x-y|^{d-\alpha}}
 \left(1\wedge \frac{  \delta_D(x)}{ |x-y| \wedge 1}\right)^{\alpha/2}  \left(1\wedge \frac{  \delta_D(y)}{ |x-y| \wedge 1 }\right)^{\alpha/2} \\
& \le G^m_D(x, y)    \le c \frac{1 +
(m^{1/\alpha} |x-y|)^{2-\alpha}     }{|x-y|^{d-\alpha}}
 \left(1\wedge \frac{  \delta_D(x)}{ |x-y| \wedge 1}\right)^{\alpha/2}  \left(1\wedge \frac{  \delta_D(y)}{ |x-y| \wedge 1 }\right)^{\alpha/2}.
\end{align*}
\end{thm}

Taking $m\downarrow 0$, the estimates in Theorem \ref{T:1.3} recover
the sharp Green function estimates for symmetric $\alpha$-stable
processes in exterior $C^{1,1}$ open sets when $d\geq 3$ that
was previously established in \cite{CT} for any dimension $d\geq 2$.

The rest of the paper is organized as follows.
In Section \ref{S:2}, we  summarize some basic properties of
relativistic stable processes and  give the proof
of the upper bound estimate in Theorem \ref{t:MAIN}.
In Section \ref{S:ub2}, we present interior lower bound estimates
for $p^m_D$ in exterior open sets.
Lower bound estimates for $p^m_D(t, x, y)$
up to the boundary are established in
Section \ref{sec:exterior1} for $t\leq T/m$ and in Section
\ref{sec:exterior2} for $t>T/m$.
The proof of Theorem \ref{T:1.3} is given
in Section \ref{S:6}.

Throughout this paper, we assume that
$\alpha\in (0, 2)$
and $m>0$. The values of the constants $C_1, C_2, C_3$ will remain the same throughout this paper,
while $c_1, c_2, \cdots$ stand for constants whose values are unimportant and which may change from location to location.
The labeling of the constants $c_1, c_2, \cdots$ starts anew in the proof of each result. The
dependence of the constant $c$ on the dimension $d$ will not be
mentioned explicitly. We will use ``$:=$" to denote a definition,
which is read as ``is defined to be". We will use
$\partial$ to denote a cemetery point and for every function $f$, we
extend its definition to $\partial$ by setting $f(\partial )=0$. We
will use $dx$ to denote the Lebesgue measure in $\bR^d$. For a Borel
set $A\subset \bR^d$, we also use $|A|$ to denote its Lebesgue
measure and $aA:=\{ ay: y \in A\}$ for $a>0$.
For two non-negative functions $f$ and $g$, the
notation $f\asymp g$ means that there are positive constants $c_1,
c_2$ so that $c_1g(x)\leq f (x)\leq c_2 g(x)$
in the common domain of definitions
 for $f$ and $g$.

\section{Basic properties of relativistic stable processes}
\label{S:2}

A symmetric $\alpha$-stable process $X=\{X_t, t\geq 0, \P_x, x\in
\bR^d\}$ in
$\bR^d$, where $d\geq 1$,
is a L\'evy process
whose characteristic
function is given by
 \eqref{e:ch} with $m=0$.
The L\'evy density of $X$ is given by
$ J(x)=j(|x|)= {\cal A}(d, -\alpha)|x|^{-(d+\alpha)}$,
where
$$
{\cal A}(d, -\alpha)=\frac{ \alpha  \, \Gamma(\frac{d+\alpha}2)}
{2^{1-\alpha} \, \pi^{d/2}\Gamma(1-\frac{\alpha}2)}.
$$
Here $\Gamma$ is the Gamma function defined by $\Gamma(\lambda):=
\int^{\infty}_0 t^{\lambda-1} e^{-t}dt$ for every $\lambda > 0$.

The L\'{e}vy measure of $X^m$ has a density
 \begin{equation}\label{e:jm2}
J^m(x)=j^m(|x|)=  {\cal A} (d, \, -\alpha) |x|^{-d-\alpha} \psi
(m^{1/\alpha}|x|)=j(|x|)\psi
(m^{1/\alpha}|x|)
\end{equation}
where
 \bee\label{e:psi} \psi (r):= 2^{-(d+\alpha)} \, \Gamma \left(
\frac{d+\alpha}{2} \right)^{-1}\, \int_0^\infty s^{\frac{d+\alpha}{
2}-1} e^{-\frac{s}{ 4} -\frac{r^2}{ s} } \, ds,
 \eee
which is decreasing and a smooth function of $r^2$ satisfying
 $\psi (0)=1$ and
 \begin{equation}\label{e:2.10}
\psi (r) \asymp \phi(r):=e^{-r}(1+r^{(d+\alpha-1)/2} ) \qquad \hbox{on } [0,
\infty)
\end{equation}
(see \cite[Lemma 2]{R} and \cite[pp. 276--277]{CS4} for the details).

Put $J^m(x,y):= j^m(|x-y|)$.
The L\'evy density gives rise to a L\'evy system for $X^m$, which
describes the jumps of the process $X^m$:
for any $x\in \R^d$, stopping time $T$ (with respect to the filtration of
$X^m$) and non-negative Borel function $f$ on $\bR_+ \times \bR^d\times \bR^d$
with $f(s, y, y)=0$ for all $y\in \R^d$,
\begin{equation}\label{e:levy}
\E_x \left[\sum_{s\le T} f(s,X^m_{s-}, X^m_s) \right]= \E_x \left[
\int_0^T \left( \int_{\bR^d} f(s,X^m_s, y) J^m(X^m_s,y) dy \right)
ds \right].
\end{equation}
(See, for example, \cite[Appendix A]{CK2}.)

We will use $p^m(t, x, y)=p^m(t, x-y)$ to denote the transition
density of $X^m$.
From \eqref{e:ch}, one can easily see that $X^m$ has the following
approximate scaling property: for every $b>0$
\begin{equation}\label{e:scaling}
\left\{ b^{-1/\alpha} \big( X^{m/b}_{b t}-X^{m/b}_0 \big),
 t\geq 0\right\}
\hbox{ has the same distribution as that of }
\left\{ X^m_t-X^m_0,t\geq 0 \right\}.
\end{equation}
In terms of transition densities,
this scaling property can be written as
\bee\label{scale_p}
p^m(t,x,y) \,=\, b^{d/\alpha}p^{m/b} (bt, b^{1/\alpha} x, b^{1/\alpha}y) \quad \text{for every }t,b>0, x,y \in \R^d.
\eee
For any $m, c>0$, we define a function
$\widetilde\Psi_{d, \alpha, m, c} ( t, x, y)$
 on $(0, \infty)\times \R^d\times \R^d$ by
$$
\widetilde\Psi_{d, \alpha, m, c} ( t, x, y)
:= \begin{cases}
 {t^{-d/\alpha}}\wedge   tJ^m(x, y)  ,
 ~ &\forall t \in
(0,1/m];
\\  m^{d/\alpha-d/2}t^{-d /2}\exp\left( -c^{-1} (m^{1/\alpha}
|x-y|\wedge m^{2/\alpha-1}\frac{|x-y|^2}{t}   )\right), ~ &\forall t \in
(1/m, \infty).
\end{cases}
$$
Using \cite[Theorem 1.2]{CKK3}, \cite[Theorem 4.1]{CKS3} and \eqref{scale_p} we get

\begin{thm}\label{T1u}
There exist $c_1, C_1>1$ such that for all $m>0$ and
$(t, x, y)\in (0, \infty)\times \R^d\times \R^d$,
$$
c_1^{-1}\widetilde\Psi_{d, \alpha, m, 1/C_1} (t, x, y)\le
p^m(t,x,y)\le c_1\widetilde\Psi_{d, \alpha, m, C_1} ( t, x, y).
$$
\end{thm}

For any open set $D$, we use $\tau^m_D:=\inf\{t>0: \, X^m_t\notin D\}$ to
denote the first exit time from $D$ by $X^m$, and $X^{m,D}$ to denote the subprocess of $X^m$
killed upon exiting $D$ (or, the killed relativistic stable process
in $D$ with mass $m$).
It is   known (see \cite{CK2}) that $X^{m,D}$
has a continuous transition density $p^m_D(t, x, y)$ with respect to
the Lebesgue measure. $p^m_D(t, x, y)$ has the following
scaling property:
 \bee\label{scale_kp}
p_D^m(t,x,y) \,=\, b^{d/\alpha}p_{b^{1/\alpha}D}^{m/b} (bt, b^{1/\alpha}
x, b^{1/\alpha} y)
\quad \text{for every }t,b>0, x,y \in D.
 \eee
Thus the Green function $G^m_D(x, y):=\int_0^\infty p^m_D(t, x, y)dt$
of $X^{m,D}$ satisfies
 \bee\label{scale_kg} G^m_{D}(x, y) =
b^{(d-\alpha)/\alpha} G^{m/b}_{b^{1/\alpha}D} (b^{1/\alpha} x,
b^{1/\alpha} y) \qquad
\text{for every } b>0, x,y \in D.
 \eee

We now introduce the space-time process $Z^m_s:=(V_s, X^m_s)$, where
$V_s=V_0- s$.
The law of the space-time process $s\mapsto Z^m_s$ starting from $(t,
x)$ will be denoted as $\mathbb{P}^{(t, x)}$ and as usual,
 $\mathbb{E}^{(t, x)}
 [ \,\cdot\,]
  =\int \,\cdot\, \mathbb{P}^{(t, x)} (d \omega).$

We say that a non-negative Borel function $h(t,x)$ on $[0,
\infty)\times \bR^d$ is {\it parabolic} with respect to the process
$X^m$ in a relatively open subset $E$ of $[0, \infty)\times \bR^d$
if for every relatively compact open subset $E_1$ of $E$,
$h(t, x)=\E^{(t,x)}
\left[ h \big( Z^m_{\widetilde \tau^{m}_{E_1}} \big) \right]$
for every $(t, x)\in E_1, $ where
$\widetilde\tau^m_{E_1} =\inf\{s> 0: \, Z^m_s\notin E_1\}$.
Note that $p_D^m(\cdot, \cdot, y)$ is parabolic with respect to the
process $X^m$ in $(0, \infty)\times D$.

The following uniform parabolic Harnack inequality is
an extension of
 \cite[Theorem 2.9]{CKS3} in that
it is stated for all $r>0$ and $m>0$ instead of only for $r\in (0, R]$ and $m\in (0, M]$.
Due to the recent result in \cite{CKK3},
the following uniform parabolic Harnack inequality  is an easy consequence of the approximate scaling
property \eqref{e:scaling}  and
the parabolic Harnack inequality \cite[Theorem 4.11]{CKK3}.

\begin{thm}\label{T:2.3}
For $M>0$ and $\delta\in (0, 1)$,  there exists
$c=c(d, \alpha, \delta, M)>0$ such that for every
$m>0$, $x_0\in \bR^d$,
$t_0\ge 0$,
$r>0$  and every non-negative function $u$ on $[0,
\infty)\times \bR^d$ that is
  parabolic with respect to the
process $X^m$ on
  $(t_0,t_0+4\delta
  (r^\alpha\vee m^{2/\alpha-1} r^2)
  ] \times B(x_0,4r)$,
$$
\sup_{(t_1,y_1)\in Q_-}u(t_1,y_1)\le
c \, \inf_{(t_2,y_2)\in
Q_+}u(t_2,y_2),
$$
where $Q_-=[t_0+
\delta
  (r^\alpha\vee m^{2/\alpha-1} r^2),t_0+2\delta
  (r^\alpha\vee m^{2/\alpha-1} r^2)]\times B(x_0,r)$
and $Q_+=[t_0+3\delta
  (r^\alpha\vee m^{2/\alpha-1} r^2),t_0+ 4\delta
  (r^\alpha\vee m^{2/\alpha-1} r^2)]\times B(x_0,r)$.
\end{thm}

We now prove the upper bound estimate in Theorem \ref{t:MAIN}.

\medskip

\noindent {\bf Proof of the upper bound estimate in Theorem \ref{t:MAIN}.}
Without loss of generality, we assume
$M=1/3$
and $T=1$.
In view of Theorems \ref{t:main}(i),  we only need to prove the upper bound in Theorem \ref{t:MAIN} 
for $t\ge 3$.
By the semigroup property and Theorem \ref{t:main}(i),
we have for $t\geq 3$, $0<m \le 1/3$
and $x, y\in D$,
\begin{eqnarray}
    p^m_D(t,x,y)
 &=& \int_D \int_D p^m_D(1,x,z)p^m_D(t-2,z,w)p^m_D(1,w,y)dzdw \nn\\
&\le & c_1\,  (1\wedge \delta_D(x))^{\alpha/2} (1\wedge
\delta_D(y))^{\alpha/2} f(t, x, y), \label{e:fgwe1}
\end{eqnarray}
where
\begin{eqnarray*}
f(t, x, y)=\int_{\R^d \times \R^d}
 \left( 1 \wedge \frac{ \phi ( m^{1/\alpha}
  |x-z|/
(16)) }{|x-z|^{d+\alpha}} \right)
  p^m(t-2,z,w) \left( 1 \wedge \frac{ \phi ( m^{1/\alpha}
  |w-y|/
(16)) }{|w-y|^{d+\alpha}} \right) dzdw.
\end{eqnarray*}
By Theorem \ref{T1u} and \eqref{e:2.10},
 there exists a constant $A\ge 16$
 such that for every $t \ge 3$
\begin{align*}
& f(t, x, y)\\
  \le & c_2\,  \int_{\R^d \times \R^d}\left( 1 \wedge \frac{ \phi ( A^{-1} m^{1/\alpha}
  |x-z|) }{|x-z|^{d+\alpha}} \right)
  p^m(t-2,A^{-1} z,A^{-1} w) \left( 1 \wedge \frac{ \phi (A^{-1}  m^{1/\alpha}
  |w-y|) }{|w-y|^{d+\alpha}} \right)dzdw\nn\\
   \le & c_3\,  \int_{\R^d \times \R^d} p^m(1,A^{-1} x,A^{-1} z)
  p^m(t-2,A^{-1} z,A^{-1} w) p^m(1,A^{-1} w,A^{-1} y) dzdw.
   \end{align*}
Thus by the change of variables $\wh z=A^{-1}z$, $\wh w=A^{-1}w$ and the semigroup property, we have that
   \begin{eqnarray}
  f (t, x, y)
  &\le & c_4\,    \int_{\R^d \times \R^d} p^m(1,A^{-1} x,\wh z)
  p^m(t-2,\wh z,\wh w) p^m(1,\wh w,A^{-1} y)d\wh z d\wh w
 \nn\\
      &= & c_4\,   p^m(t, A^{-1} x , A^{-1} y) .\label{e:question}
 \end{eqnarray}
 Now using \eqref{scale_kp} and  Theorem \ref{T1u} again, we conclude
  that for every $m \le 1/3$
    \begin{align*}
 &p^m(t, A^{-1} x , A^{-1} y)\,=\,
   A^{d} p^{mA^{-\alpha}} (A^{\alpha} t, x, y) \\
      \le& c_5
      \begin{cases}
 {t^{-d/\alpha}}\wedge  tj^m(c_6|x-y|),
 ~ &\forall t \in
[3,1/m];
\\ m^{d/\alpha-d/2} t^{-d /2}\exp\left( -c_6(m^{1/\alpha}
|x-y|\wedge m^{2/\alpha-1}\frac{|x-y|^2}t   )\right), ~
&\forall t >1/m. \end{cases}
  \end{align*}
This together
with \eqref{e:fgwe1} and \eqref{e:question} establishes
the upper bound estimate in Theorem \ref{t:MAIN}.
\qed

\section{Interior lower bound estimates}\label{S:ub2}

Throughout this section, we assume the dimension $d\geq 1$. We discuss
interior lower bound estimates for the heat kernel $p^m_D(t, x, y)$ all $t>0$.
We first establish interior lower bound estimates for the heat kernel $p^m_D(t, x, y)$
of an arbitrary open set for all $m>0$ and $t \le T/m$.

\begin{prop}\label{step31}
 Suppose that $D$ is an arbitrary open set in $\R^d
 $
 and $T>0$ is a constant.
  There exists a constant $c=c(d, \alpha,  T)>0$ such that for all
  $m>0$,
$(t, x, y)\in (0, T/m]\times D\times D$ with
$\delta_D(x) \wedge \delta_D (y)  \ge t^{1/\alpha}$,
we have
$p^m_D(t, x, y)\,\ge  \,c \, (t^{-d/\alpha} \wedge  { t}{J^m(x,y)})$.
\end{prop}

\pf By \cite[Proposition 3.5]{CKS3}, there is a constant $c=c(d, \alpha,  T)>0$ such that
for $m>0$ and
$(t, x, y)\in (0, T]\times D\times D$ with
$\delta_D(x) \wedge \delta_D (y) \ge t^{1/\alpha}$,
we have
$p^1_D(t, x, y)\,\ge  \,c \, (t^{-d/\alpha} \wedge  { t}{J^1 (x,y)})$.
The conclusion of the proposition for general $m>0$ follows immediately
from this and the scaling property
\eqref{scale_kp}.
 \qed

For notational convenience, we denote
the ball
$B(0, r)$ by $B_r$.
In the rest of this section,  we will establish interior lower bound estimate
on the heat kernel
$p^m_{\overline B_R^c }(t, x, y)$
for $m>0$, $R>0$, $t\ge T/m$, where
$T$ is a positive constant.
To achieve this, we first establish some
results for a large class of open sets
which might be of independent interest.

\begin{lemma}\label{L:4.2}
For any positive constants  $T$ and $a$, there exists
$c=c(d, \alpha, a, T)>0$ such that
for any
$t \ge T$,
$$
\inf_{y\in\bR^d}
\P_y \left( \tau^1_{B(y, a \sqrt{t})} >
 t \right)\, \ge\, c.
$$
\end{lemma}

\pf This result is an easy consequence of \cite[Theorem 4.8]{CKK3}. In fact, by \cite[Theorem 4.8]{CKK3}
\begin{eqnarray*}
\P_y \left( \tau^1_{B(y,a \sqrt{t})} >
t\right)
&=& \int_{B(y, a\sqrt{t} )}
p ^1_{B(y, a\sqrt{t} )}(t ,y,w) dw\\
&\ge& \int_{B(y, a\sqrt{t}/2 )} p ^1_{B(y, a \sqrt{t} )}(t ,y,w) dw  \ge  c.
\end{eqnarray*}
 \qed

\begin{lemma}\label{l:st3_3}
Let $a$ and $T$ be positive constants.
There exist $c_i=c_i(d, \alpha, T, a)>0$, $i=1,2,$ such that for all
$(t, x, y)\in [T, \infty)\times D\times D$ with $ \delta_D(x) \wedge \delta_D
(y) \ge a \sqrt{t}$ and $|x-y| \ge 2^{-1} \sqrt{t}$,
$$
\P_x \left( X^{1,D}_t \in B \big( y, \,  (a \wedge1)2^{-1} \sqrt{t}
\big) \right) \,\ge \,
c_1\,\frac{t^{1+d/2}e^{-c_2|x-y|}}{|x-y|^{d+\alpha} }.
$$
\end{lemma}

\pf By Lemma \ref{L:4.2}, starting at $z\in B(y, \, (a \wedge1)
\sqrt{t}/4)$, with probability at least $c_1=c_1(d, \alpha, T,
a)>0$, the process $X^1$ does not move more than $(a
\wedge1)6^{-1}\sqrt{t} $ by time $t $. Thus, it is sufficient to show
that there exist  constants $c_2=c_2(d, \alpha, T, a)>0$, $i=2, 3,$ such that
for any $t  \ge T$ and  $(x, y)\in D\times D$ with $|x-y| \ge
  \sqrt{t}/2$,
$$
\P_x \left(X^{1,D} \hbox{ hits the ball } B(y, \, (a
\wedge1) \sqrt{t}/4)\mbox{ by time } t \right)\, \ge\, c_2\,
\,\frac{ t^{1+d/2}e^{-c_3|x-y|}}{|x-y|^{d+\alpha} }.
$$

Let $B_x:=B(x, \, (a \wedge1)6^{-1} \sqrt{t})$, $B_y:=B(y, \, (a
\wedge1) 6^{-1} \sqrt{t})$ and $\tau^1_x:=\tau^1_{B_x}$. It follows
from Lemma \ref{L:4.2} that there exists
$c_4=c_4(d, \alpha, a, T)>0$
such that
\begin{equation}\label{e:newer2}
\E_x \left[ t \wedge\tau^1_{x} \right] \,\ge\,
t \, \P_x\left(\tau^1_{x} \ge  t\right)\,\ge\,
c_4\,t  \qquad \hbox{for }
t\ge T.
\end{equation}
By the L\'evy system in \eqref{e:levy},
\begin{eqnarray*}
&&\P_x \left(X^{1,D} \mbox{ hits the ball } B \big( y, \,
(a \wedge1)  \sqrt{t}/4 \big)\mbox{ by time }   t \right) \nonumber\\
&\ge & \P_x(X^1_{t\wedge \tau^1_x}\in B \big(y,\,(a \wedge1)
\sqrt{t}/4 \big) \hbox{ and }
t \wedge \tau^1_x \hbox{ is a jumping time })\nonumber \\
&\geq & \E_x \left[\int_0^{t\wedge \tau^1_x} \int_{B_y} J^1(X^1_s,
u) duds \right].\end{eqnarray*}
Since for any $(z,u) \in B_x \times B_y$ we have
$$
|z-u| \le |z-x|+|x-y| + |y-u| \le (a \wedge1)  \sqrt{t}/3 +|x-y|
\le (1+ 2(a \wedge1)/3 )
|x-y|,
$$
we get that
$$
\int_{B_y}J^1(X^1_s, u) du \ge
c_5 |B_y|
\frac {e^{-c_6|x-y|}}{|x-y|^{d+\alpha}}
\quad \text{ for every }
s < t \wedge \tau^1_x.
$$
Thus by \eqref{e:newer2},
\begin{align*}
\P_x \left(X^{1,D} \mbox{ hits the ball } B \big(y, \,
(a \wedge1) \sqrt{t}/4 \big)\mbox{ by time }   t \right) \,&\ge \,
c_7 \,\E_x \left[ t\wedge \tau_x \right]
|B_y| \, \frac{e^{-c_6|x-y|}}{
|x-y|^{d+\alpha}}\\
 &\ge \,
c_8\,\frac
{t^{1+d/2}e^{-c_6|x-y|}}{|x-y|^{d+\alpha}}.
\end{align*}
\qed

For an open set $D\subset \bR^d$ and $(\lambda_1,
\lambda_2)\in (1, \infty)\times (0, \infty)$, we say {\it the path
distance in
$D$ is comparable to the
Euclidean distance with characteristics $(\lambda_1, \lambda_2)$} if
the following holds for any $r >0$: for every $x, y$ in the same
component of $D$ with $\delta_D(x) \wedge \delta_D(y)\geq r$, there
is a length parameterized rectifiable curve $l$ in $D$ connecting
$x$ to $y$ so that the length of $l$ is no larger than
$\lambda_1|x-y|$ and $ \delta_D(l(u))\ge \lambda_2 r, u\in [0,
|l|] $.

Clearly, such a property holds for all Lipschitz domains
 with compact complements and domains above graphs of Lipschitz
functions.

\begin{prop}\label{step1}
Suppose that $D$ is a domain such that the path
distance in $D$ is comparable to the Euclidean distance with
characteristics $(\lambda_1, \lambda_2)$. For any positive constants
$a$ and $T$, there exists a positive constant
$c=c(d, \alpha, T, a, \lambda_1, \lambda_2)$ such that
 for all $(t, x, y)\in [T, \infty)\times D\times D$ with
$\delta_D(x)\wedge \delta_D(y) \ge a \sqrt{t}$ and $\sqrt{t} \geq
2|x-y|$,
we have $
p^1_D(t,x,y) \,\ge\,c\, t^{-d/2}.
$
\end{prop}

\pf Let  $t \ge T$ and $x, y \in D$ with
$\delta_D(x)\wedge \delta_D(y) \ge a \sqrt{t}$  and  $\sqrt{t} \geq 2|x-y|$.
The assumption that $D$ is a domain such that the path
distance in $D$ is comparable to the Euclidean distance with
characteristics $(\lambda_1, \lambda_2)$ enables us to
apply the parabolic Harnack inequality (Theorem \ref{T:2.3})
$N=N(a, \lambda_1, \lambda_2)$ times and to get that
there exists $c_1=c_1(d, \alpha, T, a, \lambda_1, \lambda_2)>0$ such
that
$$
p^1_D(t/2, x, w) \, \le  \, c_1\, p^1_D(t,x,y)\qquad \hbox{for } w
\in B\big(x, 2(a \wedge1) \sqrt{t}/3 \big) .
$$
This together with Lemma \ref{L:4.2} yields that
\begin{eqnarray*}
p^1_D(t, x, y) &\geq & \frac{1}{c_1 \, | B(x, (a \wedge1) \sqrt{t}/2)|}
\int_{B(x,(a \wedge1) \sqrt{t}/2)} p^1_D(t/2, x, w)dw\\
&\geq &
c_2 t^{-d/2} \, \P_x \left( \tau^1_{B(x,(a \wedge1)\sqrt{t}/2)} >
t/2\right) \, \ge\, c_3 t^{-d/2}
\end{eqnarray*}
where $c_i=c_i(d, \alpha, T,  a, \lambda_1, \lambda_2)>0$, $i=2, 3$. \qed

\begin{prop}\label{step3}
Suppose that $D$ is a domain such that the path
distance in $D$ is comparable to the Euclidean distance with
characteristics $(\lambda_1, \lambda_2)$. For any positive constants
$a$ and $T$,
there exist constants $c_i=c_i(d, \alpha, a, T, \lambda_1, \lambda_2)>0,$ $i=1,2,$ such that
\begin{equation}\label{e:lu3}
p^1_D(t, x, y)\,\ge  \,
c_1\,\frac{ t e^{-
c_2|x-y|}}{|x-y|^{d+\alpha} }
\end{equation}
for every $(t, x, y)\in
[T, \infty)\times D\times D$ with $\delta_D(x) \wedge \delta_D
(y) \ge a \sqrt{t}$ and $|x-y|^2 \ge
  t/8$.
\end{prop}

\pf By the semigroup property,  Proposition \ref{step1} and Lemma
\ref{l:st3_3}, there exist positive constants
$c_i=c_i(d, \alpha, T, a)$, $i=1,2,3,$
such that
\begin{eqnarray*}
p^1_D(t, x, y)
&\ge& \int_{B(y, \, (a \wedge1) 2^{-1} (t/2)^{1/2})} p^1_{D}(t/2, x,
z) p^1_{D}(t/2, z, y) dz\\
&\ge& c_1 t^{-d/2} \bP_x \left(
X^{1,D}
_{t/2} \in B(y, (a \wedge1)2^{-1}
(t/2)^{1/2}) \right)
 \,\ge \, c_2\,\frac{ t e^{-
 c_3|x-y|}}{|x-y|^{d+\alpha} }.
\end{eqnarray*}
 \qed

\begin{thm}\label{main50}
Suppose that $D$ is
a domain such that the path distance in $D$ is comparable to
the Euclidean distance with characteristics $(\lambda_1,
\lambda_2)$.
For any $C^*, a>0$, there exist $c_i=c_i(d, \alpha, a,
C^*,
\lambda_1, \lambda_2)>0$, $i=1,2$, such that for every $t\in
(0,\infty)$
and
$x, y\in D$
 with $\delta_D(x)
\wedge \delta_D (y) \ge a \sqrt{t}$,
$$
p^1_D(t, x,y)\geq c_1  t^{-d /2}\exp \left( -\frac{c_2 |x-y|^2}t \right)
\quad \hbox{when } C^* |x-y| \le t\le |x-y|^2 .
$$
\end{thm}
\pf  Fix $C^* >0$. Suppose that $x, y$ are
in $D$ with $\delta_D(x)
\wedge \delta_D (y) \ge a \sqrt{t}$ and satisfy
$C^* |x-y| \le t\le |x-y|^2$.
For simplicity,  let $R:=|x-y|$. Note that $t \ge (C^*)^2$.

By our assumption on $D$, there is a length parameterized curve $l\subset D$
connecting $x$ and $y$ such that the total length $|l|$ of $l$
is less than or equal to
$\lambda_1R$ and $ \delta_D(l(u))\ge \lambda_2a\sqrt{t}$ for every $u\in
[0, |l|]$.
Let $\lambda_3 \ge \max\{ 4/(\lambda_2^2 a^2),  \,  (12\lambda_1)^2\}$
and  $k$  the smallest integer satisfying
$
k\geq \lambda_3R^2/t.
$
(The integer $k$ depends on $t$ and $R$.)
 Then, since $t\in [C^*R, R^2]$,
\begin{align}
\frac{t}{k} \ge
\frac{t}{1+\lambda_3R^2/t}=\frac{t^2}{t+\lambda_3R^2} \ge
  \frac{( t/R)^2}{1+\lambda_3}   \ge
  \frac{(C^*)^2}{1+\lambda_3} .\label{e:pp1}
\end{align}
 Let $x_j=l( j|l|/k)$ and $B_j:=B(x_j,
\sqrt{t/k}/8)$, $j=0,1,...,k$.
 Note that, since $
\lambda_2^2 a^2/4 \leq \lambda_3
\le \lambda_3 R^2/t  \le k
$, we have
$\delta_D(x_j) > \lambda_2 a \sqrt{t}\geq 2 \sqrt{t/k}$ for each $j$.
So we have $B_j\subset D$ and for each $y\in B_j$,
$B(y, \sqrt{t/k}) \subset D$.

Observe that for $(y_j, y_{j+1}) \in B_j \times
B_{j+1}$, since $\lambda_3 > (12\lambda_1)^2$,
\begin{align}
&|y_j-y_{j+1}|\,\leq\,
|x_j-x_{j+1}|+|y_j-x_j|+|y_{j+1}-x_{j+1}|
\,\leq\, \frac{|l|}{k}+\frac{1}{4}\sqrt{t/k}\,\le\,
\frac{\lambda_1}{\sqrt{k}}\frac{R}{\sqrt{k}}+\frac{1}{4}\sqrt{t/k}\nn\\
 &\leq\, \frac{\lambda_1}{\sqrt{k}}\frac{R\sqrt{t}}
 {\sqrt{\lambda_3}R}+\frac{1}{4}\sqrt{t/k}
 = (\lambda_1/\sqrt{\lambda_3}+1/4) \sqrt{t/k}  \,<\, \sqrt{t/k}/3.
 \label{y_kk}
\end{align}
Now using \eqref{e:pp1}, \eqref{y_kk} and Proposition \ref{step1}, we get

\begin{equation}
p^1_D (t/k, y_j,y_{j+1})\geq c_1(t/k)^{-d/2}, \quad \text{for every
} (y_j,y_{j+1}) \in B_j \times B_{j+1}. \label{chalb22}
\end{equation}
Using \eqref{chalb22}  and the fact $k\geq \lambda_3R^2/t$, we have
\begin{align*}
&p^1_D(t, x,y)\,\geq \,\int_{B_1}\dots\int_{B_{l-1}} p^1_D (t/k,
x,y_1)\dots p^1_D( t/k, y_{k-1},y)dy_1\dots
dy_{k-1}\nn\\
&\geq \, c_1(t/k)^{-d/2} \Pi_{i=1}^{k-1} \left(
c_1 8^{-d}|B(0,1)|(t/k)^{-d/2}(t/k)^{d/2} \right)
= \, c_1(t/k)^{-d/2} \left(
c_1 8^{-d}|B(0,1)| \right)^{k-1}\\
&\ge \, c_2(t/k)^{-d/2}\exp (-c_3k)
\,\ge \, c_4t^{-d/2}\exp \left( -\frac{c_5|x-y|^2}t \right).
\end{align*}

\qed

Combining Theorem \ref{main50} with Propositions \ref{step1}
and \ref{step3}, we have the following lower bound estimates for
$p^1_D(t,x,y)$.

\begin{thm}\label{T:lbRd_1}
Let $a$ and $T$ be positive constants. Suppose that $D$ is a
domain such that the path distance in $D$ is
 comparable to the Euclidean distance with characteristics
$(\lambda_1, \lambda_2)$. Then there exist constants
$c_i=c_i(d, \alpha, a, T, \lambda_1, \lambda_2)>0$, $i=1,2,$
such that
for every $(t, x, y)\in [T, \infty)\times D\times D$ with
$\delta_D(x) \wedge \delta_D (y) \ge a \sqrt{t}$,
$$
p_D^1(t, x, y) \ge c_1t^{-d /2}\exp\left( -c_2\Big(|x-y| \wedge
\frac{|x-y|^2}t\Big)\right) .
$$
\end{thm}

Observe that any exterior ball
$\overline B^c$ is a domain in which the path distance
is comparable to the Euclidean distance with characteristics
$(\lambda_1, \lambda_2)$ independent of the radius of the ball $B$.
The following follows immediately from Theorem \ref{T:lbRd_1} and the scaling property \eqref{scale_kp}.

\begin{thm}\label{T:lbRd}
Let $a $ and $T$ be positive constants. Then there exist constants
$c_i=c_i(d, \alpha, a, T )>0$, $i=1,2,$
such that
for every $R>0$, $m >0$
and  $(t, x, y)\in [T/m, \infty)\times \overline B_R^c \times \overline B_R^c $ with
$\delta_{\overline B_R^c }(x) \wedge \delta_{\overline B_R^c } (y) \ge a\, m^{1/2 -1/\alpha}\sqrt{t}$,
$$
p_{\overline B_R^c }^m(t, x, y) \ge c_1m^{d/\alpha-d/2} t^{-d /2}\exp\left( -c_2
\big(m^{1/\alpha}
|x-y|\wedge m^{2/\alpha-1}\frac{|x-y|^2}t   \big)\right) .
$$
\end{thm}

\section{Small time lower bound estimates}
\label{sec:exterior1}

In the remainder of this paper we will always assume that the dimension $d \ge 3$.
The goal of this section is to establish the lower bound estimates in Theorem
\ref{t:MAIN} for $t\le T/m$, where $T$ is a positive constant.

Let $G^m(x, y)$ be the Green function of $X^m$.
It follows from \cite[Theorems 3.1 and 3.3]{RSV} that there exists $c=c(d, \alpha)>1$ such that
$$
c^{-1} (|x-y|^{\alpha-d} +|x-y|^{2-d})\,\le\, G^1(x, y)\, \le\, c \, (|x-y|^{\alpha-d} +|x-y|^{2-d}).
$$
Using this and \eqref{scale_kg} we get that for every $m>0$ and $x,y \in \R^d$,
\begin{equation}\label{e:g_m}
c^{-1} (|x-y|^{\alpha-d} +m^{(2-\alpha)/\alpha}|x-y|^{2-d})\,\le \, G^m(x, y) \,\le\, c\, (|x-y|^{\alpha-d}
 +m^{(2-\alpha)/\alpha}|x-y|^{2-d}).
\end{equation}

For a Borel set $A$,
we use $\sigma^m_A$ to denote the first
hitting time of $A$ by $X^m$.
Recall that we denote the ball $B(0, r)$ by
$B_r$.

\begin{lemma}\label{L:C0}
There is a constant $C_2 =
C_2(d, \alpha)>1$ such that for all $R, m>0$,
\begin{eqnarray*}
&&C_2^{-1} \frac{R^d}{R^\alpha+m^{(2-\alpha)/\alpha}R^2}
\left(|x|^{\alpha-d}+ m^{(2-\alpha)/\alpha}|x|^{2-d}\right) \\
&\le& \P_x (\sigma^m_{\overline{B}_R}<\infty) \,\leq\, C_2
\frac{R^d}{R^\alpha+m^{(2-\alpha)/\alpha}R^2}
\left(|x|^{\alpha-d}+ m^{(2-\alpha)/\alpha}|x|^{2-d} \right)
\qquad \hbox{for } |x| \ge 2R.
 \end{eqnarray*}
\end{lemma}

\pf
For $|x| \ge 2R$,
 \begin{eqnarray}
\int_{\overline B_R } G^m(x, y) dy
\asymp
\int_{\overline B_R }  (|x-y|^{\alpha-d} +m^{(2-\alpha)/\alpha}|x-y|^{2-d})dy
 \asymp R^d   (|x|^{\alpha-d}+ m^{(2-\alpha)/\alpha}|x|^{2-d}).  \label{e:ne1}
 \end{eqnarray}
 On the other hand,
 for $|z| \le R$,
\begin{eqnarray*}
\int_{\overline B_R } G^m(z, y) dy
\asymp
\int_{\overline B_R }  (|z-y|^{\alpha-d} +m^{(2-\alpha)/\alpha}|z-y|^{2-d})dy\,\asymp \, R^\alpha+m^{(2-\alpha)/\alpha} R^2.
 \end{eqnarray*}
Thus, by the strong Markov property of $X^m$, for $|x|>2R$,
\begin{eqnarray}
\int_{\overline B_R } G^m(x, y) dy = \E_x \left[\int_{\overline B_R }  G^m  \Big(X^m_{\sigma^m_{\overline B_R }}, y \Big)dy; \,
   \sigma^m_{\overline B_R}<\infty \right] \asymp (R^\alpha+m^{(2-\alpha)/\alpha}R^2)\P_x (\sigma^m_{\overline B_R }<\infty).
\label{e:ne2}
 \end{eqnarray}
 Combining \eqref{e:ne1} and \eqref{e:ne2}, we arrive at the conclusion of the lemma.
\qed

The above lemma quantifies the transience of $X^m$ when dimension
$d\geq 3$, which in particular implies that
for a compact set $K$ and
a point $x$ far away from the origin, with large
probability the process started at $x$ will never visit $K$.

\begin{lemma}\label{L:epsilon_smallm}
Suppose that $a$ and $T$ are positive constants.
There exist constants $\varepsilon = \varepsilon(d,\alpha, a, T) > 0$
and $L_1 = L_1(d,\alpha, a, T) \ge
10^{4/\alpha}$ such that the following holds:
for all $R>0$,
$m>0$,
$t \in (0, T/m ] $,
$|x|>L_1R$ and $y \in B(x,a t^{1/\alpha})\cap {\overline B_R^c }$,
    \begin{align*}
        \Pr_x\left( X^{m,{\overline B_R^c }}_t \in B(y,(t/2)^{1/\alpha}) \right)
        \ge \varepsilon.
    \end{align*}
\end{lemma}

\pf
Suppose that $t \in (0, T/m] $ and $y \in B(x,a t^{1/\alpha})\cap {\overline B_R^c }$.
It follows from Theorem \ref{T1u} that there exists $c_1=c_1(
T)>0$ such that
    \begin{align*}
       & \Pr_x\left( X^m_t \in B(y,(t/2)^{1/\alpha}) \right)\,\ge\, \inf_{w \in B(0,at^{1/\alpha})} \Pr_w\left( X^m_t\in B(0,(t/2)^{1/\alpha}) \right) \\
        &\ge c_1\inf_{w \in B(0,at^{1/\alpha})}
 \int_{B(0,(t/2)^{1/\alpha})} \left( {t^{-d/\alpha}}\wedge tJ^m(z,w)  \right) dz.
    \end{align*}
 Since for $  w \in B(0,at^{1/\alpha})$ and $ z \in B(0,t^{1/\alpha})$,
  $
 m^{1/\alpha}|z-w|   \le m^{1/\alpha}((a+1) t^{1/\alpha}) \le (a+1)
T^{1/\alpha} ,
 $
 we have in view of \eqref{e:jm2} that
    \begin{align*}
    \Pr_x\left( X^m_t \in B(y,(t/2)^{1/\alpha}) \right)
     & \ge c_1\inf_{w \in B(0,at^{1/\alpha})} \int_{B(0,(t/2)^{1/\alpha})} \left( {t^{-d/\alpha}}\wedge
    \frac{      t\psi ((a+1)
T^{1/\alpha})}{|z-w|^{d+\alpha}}
      \right)  dz  \\
  &\ge c_2 \inf_{w \in B(0,a)} \int_{B(0, (1/2)^{1/\alpha})}
     \left( 1\wedge \frac{1}{|w-z|^{d+\alpha}} \right)dz \\
  &\ge c_3(a+1)^{-d-\alpha}
     |B(0, (1/2)^{1/\alpha})|=: 2 \eps.
    \end{align*}
Thus for $x \in \mathbb R^d$ and $y \in B(x,at^{1/\alpha})$, we have
    \begin{align}
        \Pr_x\left( X^m_t \not\in B(y,(t/2)^{1/\alpha}) \right)
        = 1 - \Pr_x\left( X^m_t \in B(y,(t/2)^{1/\alpha}) \right)
        \le 1-2\varepsilon.
        \label{eq:epsilon-2_smallm}
    \end{align}

Since $d \ge 3$, we may choose $L_1 \ge
10^{4/\alpha}$ so that $C_2
(L_1^{2-d}+  L_1^{\alpha-d}) \le \varepsilon$.  By Lemma \ref{L:C0},
for all $x$ with $|x| > L_1R$ we have
    \begin{align}
       &\Pr_x\left( \tau^m_{\overline B_R^c} \le t \right)
        \le  \Pr_x\left( \sigma^m_{\overline B_R} < \infty \right)
        \le C_2  \frac{R^d}{R^\alpha+m^{(2-\alpha)/\alpha}R^2}   (|x|^{\alpha-d}+ m^{(2-\alpha)/\alpha} |x|^{2-d})\nonumber\\
  &   \le C_2 \left(
     \frac{R^\alpha}{R^\alpha+m^{(2-\alpha)/\alpha}R^2} L_1^{\alpha-d}+\frac{R^2}{m^{-(2-\alpha)/\alpha}R^\alpha+R^2}
     L_1^{2-d} \right)
               \le   C_2  (L_1^{\alpha-d}+ L_1^{2-d})  \le    \varepsilon.  \label{eq:epsilon-3_smallm}
    \end{align}
    Hence, combining \eqref{eq:epsilon-2_smallm} and \eqref{eq:epsilon-3_smallm} gives
    \begin{align*}
        \Pr_x\left( X^{m, {\overline B_R^c }}_t \in B(y,(t/2)^{1/\alpha}) \right)
        &\ge \Pr_x\left( \tau^m_{\overline B_R^c} > t \right) -
        \Pr_x\left( X^{m, {\overline B_R^c }}_t \not\in B(y,(t/2)^{1/\alpha}); \, \tau^m_{\overline B_R^c} > t \right)\\
        &\ge \Pr_x\left( \tau^m_{\overline B_R^c}>t \right) - \Pr_x\left( X^m_t \not\in B(y,(t/2)^{1/\alpha}) \right) \\
        &\ge \left( 1-\varepsilon \right) - \left( 1 - 2\varepsilon \right) = \varepsilon.
    \end{align*}
\qed

\begin{lemma}\label{L:exterior-lb-large-time_smallm}
  Let $T>0$ be a constant and $L_1=L_1(d, \alpha, 3, T)$
  be the constant in  Lemma \ref{L:epsilon_smallm}.
  There exists constant $c=c(T, d, \alpha)>0$ such that for all $m>0$,
 $R>0$, $t\in (0, T/m
 ]$,
 $x, y$ satisfying
$|x|>L_1R$, $|y|>L_1R$ with $|x-y| \le (t/6 )^{1/\alpha}$, we have
$
p^m_{\overline B_R^c}
(t,x,y) \,\ge\, c\, t^{-d/\alpha}.
$
\end{lemma}

\pf
Assume without loss of
generality that $|y| \ge |x|$. If  $\delta_{\overline B_R^c }(y) \ge (2t)^{1/\alpha}$, then
$$
\delta_{\overline B_R^c }(x)\ge \delta_{\overline B_R^c }(y)-|x-y|\ge
  ( 2^{1/\alpha}-6^{-1/\alpha})  t^{1/\alpha}> t^{1/\alpha},
$$
and hence the lemma follows immediately from
Proposition \ref{step31}.

Now assume $\delta_{\overline B_R^c }(y) < (2t)^{1/\alpha}$. Since  $t-4^{-1} \delta_{\overline B_R^c }(y)^\alpha > 4^{-1} \delta_{\overline B_R^c }(y)^\alpha$,  by the semigroup property and
Theorem \ref{T:2.3}
 we have
    \begin{align}
         p^m_{\overline B_R^c }(t,x,y)
        &\ge \int_{B(y,8^{-1/\alpha} \delta_{\overline B_R^c }(y))}p^m_{\overline B_R^c }
        \left( 4^{-1} \delta_{\overline B_R^c }(y)^\alpha, x, z \right)p^m_{\overline B_R^c }\left( t - 4^{-1} \delta_{\overline B_R^c }(y)^\alpha,z,y \right)dz \notag \\
        &\ge c_1 \Pr_x\left( X^{m, \overline B_R^c }_{4^{-1} \delta_{\overline B_R^c }(y)^\alpha}
        \in B(y,8^{-1/\alpha} \delta_{\overline B_R^c }(y)) \right)p^m_{\overline B_R^c }(t-\frac38 \delta_{\overline B_R^c }(y)^\alpha,y,y).
        \label{E:exterior-lb_smallm1}
    \end{align}
Observe that
$        |x-y| \le 2|y| = 2(\delta_{\overline B_R^c }(y) + R) < 3\delta_{\overline B_R^c }(y),
 $
where the last inequality follows because $|y| > L_1R \ge 3R$
implies $\delta_{\overline B_R^c }(y) > 2R$.  Thus, $y \in B(x,3\delta_{\overline B_R^c }(y))
\cap \overline B_R^c$
and
Lemma \ref{L:epsilon_smallm} gives
$$
\Pr_x \left( X^{m, \overline B_R^c }_{4^{-1} \delta_{\overline B_R^c }(y)^\alpha}
\in B(y,4^{-1/\alpha} \delta_{\overline B_R^c }(y))\right) \ge \varepsilon.
$$

To bound the second term in \eqref{E:exterior-lb_smallm1}, we let
$s:=t -3 \cdot 8^{-1} \delta_{\overline B_R^c }(y)^\alpha$
and note that $s < t \le T/m$.
Thus, by the semigroup property, the Cauchy-Schwarz
inequality and Lemma \ref{L:epsilon_smallm},
\begin{align*}
&p^m_{\overline B_R^c }(s, y, y)\,\ge\, \int_{B(y, (s/4)^{1/\alpha})}\left(p^m_{\overline B_R^c }(s/2, y, z)\right)^2 dz\\
&\ge\,\frac1{|B(y, (s/4)^{1/\alpha})|}\left(\P_y(X^{m, D}_{s/2}\in B(y, (s/4)^{1/\alpha}))\right)^2
\,\ge\, c_2 s^{-d/\alpha}\,\ge\, c_2 t^{-d/\alpha}.
\end{align*}
The proof is now complete.
\qed

\begin{prop}\label{P:exterior-compare_smallm}
Let $T>0$ be a constant and  $L_1=L_1(d, \alpha, 3, T)$
be the constant in  Lemma \ref{L:epsilon_smallm}.
There is a constant $c=c(d, \alpha, T)>0$ such that
for every $m>0$, $R>0$, $t \in (0, T/m]$
and $x, y$ with
$|x|>L_1R$ and $|y|>L_1R$,
$$
p^m_{\overline B_R^c }(t, x, y) \geq c \left({t^{-d/\alpha}}\wedge  tj^m(2|x-y|) \right).
$$
\end{prop}

\pf  Let  $|x|>L_1R$, $|y|>L_1R$ and $t\in (0, T/m]$.
By Lemma \ref{L:exterior-lb-large-time_smallm},
we only need to show that
\begin{equation}\label{e:newcl1}
p^m_{\overline B_R^c }(t, x, y) \geq c_1 \left({t^{-d/\alpha}}\wedge  tj^m(2|x-y|) \right)
\qquad \hbox{when } |x-y|> (t/6)^{1/\alpha} .
\end{equation}
If $ t < 60\cdot 4^\alpha  R^\alpha$,
$$
\delta_{\overline B_R^c }(x) \wedge \delta_{\overline B_R^c }(y)  \ge (L_1-1)R \ge (10^{4/\alpha}-1)R
 \ge 60^{1/\alpha} \cdot 4 R > t^{1/\alpha}.
 $$
 Thus, by Proposition \ref{step31}, \eqref{e:newcl1} is true in this case.

Now we  assume that
$ t\geq 60\cdot 4^\alpha  R^\alpha$.
Since one of $|x|$ and $|y|$ should be no less
than $|x-y|/2$, we assume without loss of generality that
$|y|\ge|x-y|/2$. Let $x_0:=x+(t/60)^{1/\alpha}x /|x|$.
Note that $ B(x_0, (t/60)^{1/\alpha}) \subset \overline B_R^c $.
Since $|x-y| > ( t/6)^{1/\alpha}$, we get for every $z\in
B(x_0, (t/60)^{1/\alpha}/4)$,
$|x-z| \le|x_0-z| + (t/60)^{1/\alpha} < (t/12)^{1/\alpha}$
and
$|z-y| \le  |x-y| +|x-z| < |x-y| + (t/12)^{1/\alpha}<2 |x-y|$.
Moreover,
$$
\delta_{\overline B_R^c }(y) = |y| - R \ge \frac12|x-y|-\frac14(t/60)^{1/\alpha}
\ge \frac12(t/6)^{1/\alpha}-\frac14 (t/60)^{1/\alpha}
 >  \frac14(t/6)^{1/\alpha}
$$
 while
for $z \in B(x_0, \frac14(t/60)^{1/\alpha})$, we have $|z|\ge |x|>L_1R$ and
\begin{eqnarray*}
\delta_{\overline B_R^c }(z)
&=& |z| - R
\ge |x_0|- |x_0-z|-\frac14 (t/60)^{1/\alpha} \\
& \ge&  |x| +(t/60)^{1/\alpha}
 -\frac14 (t/60)^{1/\alpha}-\frac14 (t/60)^{1/\alpha}
 \ge \frac12 (t/60)^{1/\alpha}.
\end{eqnarray*}
Let $a=4^{-\alpha} (60)^{-1}$.
By the semigroup property,  Proposition \ref{step31} and Lemma
\ref{L:exterior-lb-large-time_smallm}, there exist positive constants
$c_i=c_i (d, \alpha, T)$ for $i=2,3$, such that
\begin{eqnarray*}
p^m_{\overline B_R^c }(t, x, y) &=& \int_{\overline B_R^c } p^m_{\overline B_R^c }((1-a)t, x, z)
p^m_{\overline B_R^c }(at, z, y)dz\\
&\ge& \int_{B(x_0, \frac14(t/60)^{1/\alpha})} p^m_{\overline B_R^c }((1-a) t, x,
z) p^m_{\overline B_R^c }(a t, z, y) dz\\
&\ge&c_2 \int_{B(x_0, \frac14(t/60)^{1/\alpha})}  t^{-d/\alpha}  \left({t^{-d/\alpha}}\wedge  tJ^m(z,y) \right) dz\\
&\ge & c_3\left({t^{-d/\alpha}}\wedge  t
j^m (2|x-y|) \right).
\end{eqnarray*}
This proves \eqref{e:newcl1}.
 \qed

We recall the following lemma from \cite{CKS5}.
\begin{lemma}[{\cite[Lemma 2.2]{CKS5}}]\label{L:ratio}
Let $\lambda, T, M$  be fixed positive
constants.
Suppose $x,x_0 \in \R^d$ satisfy $|x-x_0| = \lambda
T^{1/\alpha}$.  Then for all $a\in (0, M]$ and $z\in \R^d$,
\begin{align}
  T^{-d/\alpha} \wedge \frac{T \phi(a |x-z|)}
 {|x-z|^{d+\alpha}}
\, \asymp \,  T^{-d/\alpha} \wedge \frac{T
  \phi(a |x_0-z|)} {|x_0-z|^{d+\alpha}} ,
\label{E:ratio}
\end{align}
where the (implicit) comparison constants in  \eqref{E:ratio} depend
only on $d$, $\alpha$, $M$,  $\lambda$ and $T$.
 \end{lemma}

\begin{thm}\label{T:exterior_1}
Let $T, M$ and $R$ be positive constants.
Suppose that $D$ is an exterior $C^{1,1}$ open set in
$\R^d$ with $C^{1,1}$ characteristics $(r_0, \Lambda_0)$ and $D^c\subset B(0, R)$.
 Then there is a positive constant
$c=c(d, \alpha, r_0, \Lambda_0, R, M, T)$ so that
for all $0<m \le M$ and $(t, x, y)\in (0, T/m]\times D \times D$,
$$
p^m_D(t,x,y)\ge c
\left( 1\wedge \frac{\delta_D(x)}{1\wedge t^{1/\alpha}}\right)^{\alpha/2}
        \left( 1\wedge \frac{\delta_D(y)}{1\wedge t^{1/\alpha}} \right)^{\alpha/2}
       \left( {t^{-d/\alpha}}\wedge tj(4m^{1/\alpha}|x-y|)
       \right).
       $$
\end{thm}

\pf Without loss of generality, we assume
$M=1/3$ and $T=1$.
By Theorem \ref{t:main}, we only need to show the theorem for $t >3$.
For $x$ and $y$ in $D$, let $v \in \R^d$ be any unit
vector satisfying $
x \cdot v \ge 0$ and $
y \cdot v \ge 0$.
Let $L_1=L_1(d, \alpha, 3, 1)$ be the constant given by Lemma  \ref{L:epsilon_smallm} and define
\begin{align}\label{E:outward-push}
    x_0 := x + L_1^2Rv
    \quad\text{and}\quad
    y_0 := y + L_1^2Rv.
\end{align}
Then
$
    |x_0|^2 = |x|^2+(L_1^2R)^2+2L_1^2 \, R
\, x \cdot v
    \ge (L_1^2R)^2,
$
and similarly, $|y_0|^2 \ge (L_1^2R)^2$.

Using the semigroup property and Theorem \ref{t:main}(i), for every
 $m\in (0, 1/3]$ and  $t \in [3, 1/m]$ we have
\begin{eqnarray}\label{e:lowerf1}
   p^m_D(t,x,y)
 &=& \int_D \int_D p^m_D(1,x,z)p^m_D(t-2,z,w)p^m_D(1,w,y)dzdw \nn\\
&\ge& c_1\,  (1\wedge \delta_D(x))^{\alpha/2} (1\wedge
\delta_D(y))^{\alpha/2}
  f_1(t, x, y),
\end{eqnarray}
where
\begin{eqnarray*}
f_1(t, x, y)&=&\int_{D \times D }
  (1\wedge \delta_D(z))^{\alpha/2}
 \left( 1 \wedge \frac{ \phi ( m^{1/\alpha}
  |x-z|) }{|x-z|^{d+\alpha}} \right)  p^m_D(t-2,z,w)
 (1\wedge \delta_D(w))^{\alpha/2} \\
&&  \hskip 0.2truein \cdot \left( 1 \wedge \frac{ \phi ( m^{1/\alpha}
  |w-y|) }{|w-y|^{d+\alpha}} \right) dzdw .
\end{eqnarray*}

Note that by
 Proposition \ref{P:exterior-compare_smallm}, for $z, w\in B(0, L_1R)^c$ and $t \in [3, 1/m]$,
$$
 p^m_D(t-2,z,w) \ge  p^m_{\overline B_R^c }(t-2,z,w) \ge c_2
 \left({(t-2)^{-d/\alpha}}\wedge  (t-2)j^m(2|z-w|) \right).
$$
Since $m\leq
3$,
the lower bound
estimate in Theorem \ref{t:main}(i)
 and the above display together with
Lemma \ref{L:ratio}
 imply that
\begin{eqnarray*}
&& f_1(t, x, y)\\
  &\ge & c_3\,  \int_{D \times D }(1\wedge \delta_D(z))^{\alpha/2}
 \left( 1 \wedge \frac{ \phi ( m^{1/\alpha}
 |x_0-z|) }{|x_0-z|^{d+\alpha}} \right)
   p^m_{D}(t-2,z,w)(1\wedge \delta_D(w))^{\alpha/2}\\
&&  \hskip 0.2truein \cdot  \left( 1 \wedge \frac{ \phi ( m^{1/\alpha}
 |w-y_0|) }{|w-y_0|^{d+\alpha}} \right) dzdw\nn\\
 &\ge & c_4\,  \int_{B(0, L_1R)^c \times B(0, L_1R)^c }(1\wedge \delta_D(z))^{\alpha/2}(1\wedge \delta_D(w))^{\alpha/2}
 \left( 1 \wedge \frac{ \phi ( m^{1/\alpha}
|x_0-z|) }{|x_0-z|^{d+\alpha}} \right)
   \\
&&  \hskip 0.2truein \cdot  \left({(t-2)^{-d/\alpha}}\wedge  (t-2)j^m(2|z-w|) \right)\left( 1 \wedge \frac{ \phi ( m^{1/\alpha}
|w-y_0|) }{|w-y_0|^{d+\alpha}} \right) dzdw\nn\\
  &\ge & c_5\,  \int_{B(0, L_1R)^c \times B(0, L_1R)^c }
 \left( 1 \wedge \frac{ \phi ( 2m^{1/\alpha}
 |x_0-z|) }{|x_0-z|^{d+\alpha}} \right)
   \left({(t-2)^{-d/\alpha}}\wedge  (t-2)j^m(2|z-w|) \right)\\
&&  \hskip 0.2truein \cdot  \left( 1 \wedge \frac{ \phi (2 m^{1/\alpha}
  |w-y_0|) }{|w-y_0|^{d+\alpha}} \right) dzdw.
   \end{eqnarray*}
   Thus by the change of variables $\wh z=2z$, $\wh w=2w$,
and  Theorem \ref{T1u}, we have that
   \begin{eqnarray}
  f_1 (t, x, y)  &\ge & c_6\,    \int_{B(0, 2L_1R)^c \times B(0, 2L_1R)^c } \left( 1 \wedge \frac{ \phi ( m^{1/\alpha}
  |2x_0-\wh z|) }{|2 x_0-\wh z|^{d+\alpha}} \right) p^m(t-2,\wh z,\wh w)
 \nn\\
&&  \hskip 0.2truein \cdot
\left( 1 \wedge \frac{ \phi ( m^{1/\alpha}
  |\wh w-2y_0|) }{|\wh w-2y_0|^{d+\alpha}} \right)
  d\wh z d\wh w\nn\\
      &\ge & c_7\, \int_{B(0, 2L_1R)^c \times B(0, 2L_1R)^c }
       p^m_{B(0, 2L_1R)^c}(1, 2 x_0,    \wh z)
  p^m_{B(0, 2L_1R)^c}(t-2, \wh z, \wh w) \nn\\
&&  \hskip 0.2truein \cdot
p^m_{B(0, 2L_1R)^c}(1, \wh w, 2 y_0)
 dzdw\nn\\
      &= & c_7\,   p^m_{B(0, 2L_1R)^c}(t,2 x_0, 2 y_0).\nn
 \end{eqnarray}
We conclude from   \eqref{scale_kp} and  Proposition \ref{P:exterior-compare_smallm}  that
(recall that $|x_0|, |y_0| \ge L_1(L_1R)$)
$$
 p^m_{B(0, 2L_1R)^c}(t,2 x_0, 2 y_0)=
2^{-d}p^{m
2^{\alpha}}_{B(0, L_1R)^c}(
2^{-\alpha}
 t, x_0, y_0)
         \ge c_8 \left({t^{-d/\alpha}}\wedge  tj(4m^{1/\alpha}|x-y|) \right).
$$
Combining the last two displays with \eqref{e:lowerf1} completes the proof.
\qed

\section{Large time lower bound estimates}
\label{sec:exterior2}

The goal of this section is to establish the lower bound estimates in Theorem
\ref{t:MAIN} for $t\ge T/m$, where $T$ is a positive constant.
For any $y \in \R^d\setminus\{0\}$ and any $r>0$, we define
$$
H(y, r):=\{z\in B(y, r): z\cdot y\ge 0\}.
$$

\begin{lemma}\label{L:epsilon}
Let $T>0$. There exist constants $\varepsilon = \varepsilon(d,\alpha, T) > 0$,
$L_2=L_2(d,\alpha, T) \ge 3$ such that the following holds: for all $t \ge T$,
$R>0$, $x$ and $y$ satisfying $|x|>L_2R$, $|y|>R$ and $y \in B(x,9\sqrt{t})$,
\begin{align*}
\Pr_x\left( X^{1,\overline B_R^c }_t \in H(y,   \sqrt{t}/2) \right) \ge \varepsilon.
\end{align*}
\end{lemma}

\pf
It follows from Theorem \ref{T1u} that there exist
constants $c_1>0$ and $c_2>1$
such that
\begin{eqnarray*}
\Pr_x\left( X^1_t \in H(y,  \sqrt{t}/2 ) \right)
&\ge& \inf_{w \in B(y,9\sqrt{t})} \Pr_w\left( X^1_t\in H(y, \sqrt{t}/2) \right)\\
&\ge&
c_1\inf_{w \in B(y,9\sqrt{t})}\int_{H(y,\frac12 \sqrt{t})}
\widetilde\Psi_{d, \alpha, 1,
c_2^{-1}}( t,w, z)dz.
\end{eqnarray*}
    If $T<1$ and $t \in [T,1]$, then clearly
    \begin{align*}
     &\inf_{w \in B(y,9\sqrt{t})} \int_{H(y,\frac12\sqrt{t})}
     \widetilde\Psi_{d, \alpha, 1, c_2^{-1}}( t, w, z)dz \\
     &\ge c_3 \inf_{w \in B(y,9)} \int_{H(y,\frac{T}2)}
     \left( 1\wedge \frac{1}{|w-z|^{d+\alpha}} \right)dz \ge c_410^{-d-\alpha}|B(0,1/2)|.
    \end{align*}
     If $t >1,$ then
     \begin{align*}
          & \inf_{w \in B(y,9\sqrt{t})} \int_{H(y,\frac12\sqrt{t})}\widetilde
          \Psi_{d, \alpha, 1,
          c_2^{-1}}( t, w, z)dz \\
    & \ge \inf_{w \in B(y,9\sqrt{t})} \int_{H(y,\frac12\sqrt{t})} t^{-d /2}\exp\left( -c_2 \frac{|z-w|^2}t \right)  dz  \\
     &\ge c_5 \inf_{w \in B(y,9)} \int_{H(y,\frac12)}
     \exp\left( -c_2 {|z-w|^2} \right)dz \ge c_6e^{-c_210^2}|B(0,1/2)|.
    \end{align*}
Hence there is  $\varepsilon \in (0, 1/4)$  so that for any $t\geq T$,
$x \in \mathbb R^d$ and $y \in B(x,9\sqrt{t})$,
\begin{align}
\varepsilon<\frac12 \Pr_x\left( X^1_t \in H(y, \sqrt{t}/2) \right).
\label{eq:epsilon-2}
\end{align}

Since $d \ge 3$, we may choose $L_2 \ge 3$ so that $C_2
(L_2^{2-d}+ L_2^{\alpha-d}) \le \varepsilon$.  By Lemma \ref{L:C0},
for all $x$ with $|x| > L_2R$, we have
       \begin{align}
       &\Pr_x\left( \tau^1_{\overline B_R^c } \le t \right)
        \le \Pr_x\left( T^1_{B(0,R)} < \infty \right)
        \le C_2  \frac{R^d}{R^2+R^\alpha}   (|x|^{2-d}+ |x|^{\alpha-d}) \nonumber\\
  &   \le C_2 \left(
     \frac{R^2}{R^2+R^\alpha}   L_2^{2-d}+\frac{R^\alpha}{R^2+R^\alpha}  L_2^{\alpha-d} \right)
               \le   C_2  (  L_2^{2-d}+ L_2^{\alpha-d})  \le    \varepsilon.
        \label{eq:epsilon-3}
    \end{align}
Combining
\eqref{eq:epsilon-2} and \eqref{eq:epsilon-3} gives
\begin{align*}
    \Pr_x\left( X^{1, \overline B_R^c }_t \in H(y, \sqrt{t}/2) \right)
    &= \Pr_x\left( \tau^1_{\overline B_R^c } > t \right) -
    \Pr_x\left( X^{1, \overline B_R^c }_t \not\in H(y,\sqrt{t}/2); \, \tau^1_{\overline B_R^c } > t \right)\\
    &\ge \Pr_x\left( \tau^1_{\overline B_R^c }>t \right) - \Pr_x\left( X^1_t \not\in H(y, \sqrt{t}/2) \right) \\
    &\ge \left( 1-\varepsilon \right) - \left( 1 - 2\varepsilon \right) = \varepsilon.
\end{align*}
This proves the lemma. \qed

\begin{lemma}\label{L:exterior-lb-large-time}
Let $T>0$ and
$L_2=L_2(d, \alpha, T/8)$ be the constant in  Lemma \ref{L:epsilon}.
There exists
a constant $c=c(\alpha, d, T)>0$ such that for all $m>0$,
$R>0$,
$t \ge T/m$ and $x, y$ satisfying $|x|>L_2R$, $|y|>L_2R$, $|x-y| \le  m^{1/2-1/\alpha} \sqrt{t}/6$,
we have
\begin{align*}
p^m_{\overline B_R^c }(t,x,y) \,\ge\, c\, m^{d/\alpha-d/2} \,t^{-d/2}.
\end{align*}
\end{lemma}

\pf
We first prove the lemma for $m=1$.
Assume without loss of
generality that $|y| \ge |x|$. If  $\delta_{\overline B_R^c }(y) >
\sqrt{t}/2$, then $\delta_{\overline B_R^c }(x)\ge \delta_{\overline B_R^c }(y)-|x-y|\ge \sqrt{t}/3$, and hence the lemma follows immediately from
Theorem \ref{T:lbRd}.

Now assume $\delta_{\overline B_R^c }(y) \le  \sqrt{t}/2$.
By the semigroup property and Theorem \ref{T:2.3} we have
    \begin{align}
        p^1_{\overline B_R^c }(t,x,y)
        &\ge \int_{H(y, (t/2)^{1/2})}p^1_{\overline B_R^c }\left( t/2, x, z \right)p^1_{\overline B_R^c }\left( t/2,z,y \right)dz \notag \\
        &\ge c_1\Pr_x\left( X^{1, \overline B_R^c }_{t/2}\in H(y,(t/2)^{1/2}) \right)p^1_{\overline B_R^c }\left((t/2)- \delta_{\overline B_R^c }(y)^2,y,y \right).
        \label{E:exterior-lb}
    \end{align}
By Lemma \ref{L:epsilon} we have $\Pr_x(X^{1, \overline B_R^c }_{t/2} \in H(y,(t/2)^{1/2}))
\ge \varepsilon$.

Note that
$
     t  \ge  s:=(t/2) -  \delta_{\overline B_R^c }(y)^2
\geq t/4 \ge T/4.
 $
Hence by the semigroup property, the Cauchy-Schwarz
inequality and Lemma \ref{L:epsilon}
\begin{eqnarray*}
p^1_{\overline B_R^c }(s, y, y)&\ge& \int_{H(y, \sqrt{s}/2 )}
\big(p^1_{\overline B_R^c }(s/2, y, z)\big)^2dz\\
&\ge&\frac2{|B(y, \sqrt{s}/2  )|}\, \P_y \left(X^{1, \overline B_R^c }_{s/2}\in
 H(y, \sqrt{s}/2) \right)^2
  \,\ge \, c_2 s^{-d/2}
   \,\ge \, c_2 t^{-d/2}.
\end{eqnarray*}
Thus by \eqref{E:exterior-lb} we have
$
        p^1_{\overline B_R^c }(t,x,y) \ge c_3 t^{-d/2}.
 $

Now we consider the general case
$m>0$,
$|x-y| \le   m^{1/2-1/\alpha} \sqrt{t}/
6$,
and $t \ge T/m$. We apply \eqref{scale_kp} to the previous case and get
$$
      p^m_{\overline B_R^c }(t,x,y) =  m^{d/\alpha}  p^1_{B(0, m^{1/\alpha} R)^c}(mt,m^{1/\alpha}x,m^{1/\alpha}y)
      \ge c_3 m^{d/\alpha-d/2}  t^{-d/2}.$$
\qed

\begin{prop}\label{P:exterior-compare}
Let $T>0$ and
$L_2=L_2(d, \alpha, T/
(16))$ be the constant in  Lemma \ref{L:epsilon}.
There exist constants $c_1=c_1(\alpha, d, T)>0$ and $C_3=C_3(\alpha, d, T)>0$ such that for all
$R, m>0$, $t \ge T/m$, and $x, y$ satisfying $|x|>L_2R$, $|y|>L_2R$,
$$
p^m_{\overline B_R^c }(t, x, y) \geq c_1m^{d/\alpha-d/2} t^{-d /2}\exp\left( -C_3 \big(m^{1/\alpha}
|x-y|\wedge m^{2/\alpha-1}\frac{|x-y|^2}t   \big)\right).
$$
\end{prop}

\pf
By Lemma \ref{L:exterior-lb-large-time}, we only need to prove the proposition
for $|x-y|> \frac16 m^{1/2-1/\alpha} \sqrt{t}$, which we will assume throughout the proof.

We first prove the lemma for $m=1$.
 If  $t< (60R)^{2}$, then
$\delta_{\overline B_R^c }(x)\ge (L_2-1)R \ge 2R  > (30)^{-1} t^{1/2}$. Thus, in this case,
 the lemma follows immediately from
Theorem \ref{T:lbRd}.

 Suppose
$t\geq  T \vee (60R)^{2}$ and $|x-y| \le \frac16  \sqrt{t}$.
As one of $|x|$ and $|y|$ must be no less than
$|x-y|/2$, we assume without loss of generality that
$|y|\ge|x-y|/2$.
Let $x_0:=x+20^{-1} \sqrt{t}x /|x|$ and observe that $ B(x_0,20^{-1} \sqrt{t}) \subset \overline B^c_{|x|}\subset \overline B_R^c $.

Since $|x-y| > \frac16  \sqrt{t}$, we get for every $z\in
B(x_0,20^{-1} \sqrt{t})$.
$$
|x-z| \le|x_0-z| + \frac{1}{20} \sqrt{t} \le  \frac{1}{10}\sqrt{t}
\le \frac{1}{6}\sqrt{t/2}
$$
and
$$
|z-y| \le  |x-y| +|x-z|  \le |x-y| + \frac{1}{10}\sqrt{t} \le
2 |x-y|.
$$
Moreover, since $R <\frac1{60}\sqrt{t} $,
\begin{align*}
\delta_{\overline B_R^c }(y) = |y| - R \ge \frac12|x-y|-\frac1{60}\sqrt{t}   \ge \frac1{12}\sqrt{t}-\frac1{60}\sqrt{t} =\frac1{15}\sqrt{t}
\end{align*}
and, for $z \in B(x_0, \frac1{60}\sqrt{t})$,
\begin{align*}
\delta_{\overline B_R^c }(z) = |z| - R \ge |x_0|- |x_0-z|-\frac1{60}\sqrt{t} \ge
|x| +\frac{1}{20} \sqrt{t}
 -\frac1{60}\sqrt{t}-\frac1{60}\sqrt{t} \ge \frac1{60}\sqrt{t}.
\end{align*}
Thus by the semigroup property,
 Theorem \ref{T:lbRd}
and Lemma
\ref{L:exterior-lb-large-time}, there exist positive constants
$c_i=c_i(d, \alpha, T)$, $i=1,\dots
3$, such that
\begin{eqnarray*}
p^1_{\overline B_R^c }(t, x, y) &=& \int_{\overline B_R^c } p^1_{\overline B_R^c }(t/2, x, z)
p^1_{\overline B_R^c }(t/2, z, y)dz\\
&\ge& \int_{
B(x_0, \frac1{60}\sqrt{t})} p^1_{\overline B_R^c }(t/2, x,
z) p^1_{\overline B_R^c }(t/2, z, y) dz\\
&\ge&c_1 \int_{B(x_0, \frac1{60}\sqrt{t})}   (t/2)^{-d/2}
(t/2)^{-d /2}\exp\left( -c_2 \big(
|z-y|\wedge \frac{|z-y|^2}{t/2}   \big)\right) dz\\
&\ge & c_3t^{-d /2}\exp\left( -4c_2 \big(
|x-y|\wedge \frac{|x-y|^2}{t}   \big)\right).
\end{eqnarray*}

Now we consider the general case
$m>0$
and $t \ge T/m$. We apply \eqref{scale_kp} to the previous case and get
\begin{align*}
      &p^m_{\overline B_R^c }(t,x,y) =  m^{d/\alpha}  p^1_{B(0, m^{1/\alpha} R)^c}(mt,m^{1/\alpha}x,m^{1/\alpha}y)
    \\  &\ge c_3 m^{d/\alpha-d/2}  t^{-d/2}\exp\left( -C_3 \big(m^{1/\alpha}
|x-y|\wedge m^{2/\alpha-1}\frac{|x-y|^2}t   \big)\right).
\end{align*}
\qed

\begin{thm}\label{T:exterior}
Suppose that $M, T, R$ are positive constant, and that $D$ is an exterior $C^{1,1}$ open set
in $\R^d$ with $C^{1,1}$  characteristics $(r_0, \Lambda_0)$ and $D^c\subset B(0, R)$.
There
exist positive constants $c_i=c_i(d, \alpha, r_0, \Lambda_0, R, M,
T)$, $i=1, 2,$ such that
for all $0<m \le M$ and $(t, x, y)\in [T/m, \infty)\times D \times D$,
\begin{align*}
 &p^m_D(t,x,y)\\
&\ge c_1 \left( 1\wedge \frac{{\delta_D(x)}}{1\wedge t^{1/\alpha}} \right)^{\alpha/2} \left( 1\wedge \frac{\delta_D(y)}
{1\wedge t^{1/\alpha}} \right)^{\alpha/2}
             m^{d/\alpha-d/2}t^{-d /2}\exp\left( -4c_2 (m^{1/\alpha}
|x-y|\wedge m^{2/\alpha-1}\frac{|x-y|^2}{t}   )\right).
\end{align*}
\end{thm}
\pf
Without loss of generality, we assume $M=T=1$. By Theorem \ref{T:exterior_1}, we may assume $t> 3/m$.
For $x$ and $y$ in $D$, let $v \in \R^d$ be any unit
vector satisfying $
x \cdot v \ge 0$ and $
y \cdot v \ge 0$.
Recall that $C_1$ is the constant in Theorem \ref{T1u} and that $C_3=C_3(\alpha, d, 1)$  is the constant in
 Proposition \ref{P:exterior-compare}.
Let $A:=4 \vee (C_1C_3)$ and
$L_3:=L_2(d, \alpha, (16)^{-1})  \vee L_2(d, \alpha, (16)^{-1}A^{-\alpha})$,
where $L_2$ is given by Proposition \ref{P:exterior-compare}. Define
\begin{align*}
x_0 := x +m^{-1/\alpha} L_3^2Rv
    \quad\text{and}\quad
y_0 := y +m^{-1/\alpha} L_3^2Rv.
\end{align*}
Then
\begin{align*}
    |x_0|^2 = |x|^2+m^{-2/\alpha}(L_3^2R)^2+2m^{-1/\alpha} L_3^2R
\, x \cdot v
    \ge (m^{-1/\alpha} L_3^2R)^2 \ge (L_3^2R)^2,
\end{align*}
and similarly, $|y_0|^2 \ge (L_3^2R)^2$.
By the
semigroup property
and Theorem \ref{T:exterior_1},
for every $t \in (3/m, \infty)$,
\begin{eqnarray}\label{e:lowerf12}
   p^m_D(t,x,y)
 &=& \int_D \int_D p^m_D(1/m,x,z)p^m_D(t-2/m,z,w)p^m_D(1/m,w,y)dzdw \nn\\
&\ge& c_1\,  (1\wedge \delta_D(x))^{\alpha/2} (1\wedge
\delta_D(y))^{\alpha/2}
  f_1(t, x, y),
\end{eqnarray}
where
\begin{eqnarray*}
f_1(t, x, y)&=&\int_{D \times D }
  (1\wedge \delta_D(z))^{\alpha/2}
 \left( m^{d/\alpha} \wedge  m^{-1}\frac{ \phi ( 4m^{1/\alpha}
  |x-z|) }{|x-z|^{d+\alpha}} \right)  p^m_D(t-2,z,w)
  \\
&&  \hskip 0.2truein \cdot (1\wedge \delta_D(w))^{\alpha/2}\left( m^{d/\alpha} \wedge m^{-1}\frac{ \phi ( 4m^{1/\alpha}
  |w-y|) }{|w-y|^{d+\alpha}} \right) dzdw .
\end{eqnarray*}

By Proposition \ref{P:exterior-compare}, for $z, w\in B(0, L_3R)^c$ and
$t \in (3/m, \infty)$ (note that $t-2 \ge 1/m$),
\begin{eqnarray}
&& p^m_D(t-2/m,z,w) \ge  p^m_{\overline B_R^c }(t-2/m,z,w) \nn \\
&\ge& c_2
 m^{d/\alpha-d/2} (t-2/m)^{-d /2}\exp\left( -C_3 \big(m^{1/\alpha}
|z-w|\wedge m^{2/\alpha-1}\frac{|z-w|^2}{t-2/m}   \big)\right).\label{e:newqw0}
\end{eqnarray}
Moreover, since
$$  m^{d/\alpha} \wedge  m^{-1}\frac{ \phi ( 4m^{1/\alpha}
  r) }{r^{d+\alpha}}   =m^{d/\alpha}\left(  1 \wedge  \frac{ \phi ( 4m^{1/\alpha}
  r) }{(m^{1/\alpha}r)^{d+\alpha}} \right),
$$
we have by Lemma \ref{L:ratio},
\begin{equation}\label{e:newqw1}
  m^{d/\alpha} \wedge  m^{-1}\frac{ \phi ( 4m^{1/\alpha}
  |x-z|) }{|x-z|^{d+\alpha}} \, \ge c_3
\left(  m^{d/\alpha} \wedge  m^{-1}\frac{ \phi ( 4m^{1/\alpha}
  |x_0-z|) }{|x_0-z|^{d+\alpha}} \right)
\end{equation}
and
\begin{equation}\label{e:newqw2}
  m^{d/\alpha} \wedge m^{-1}\frac{ \phi ( 4m^{1/\alpha}
  |w-y|) }{|w-y|^{d+\alpha}} \,   \ge c_3 \left(
  m^{d/\alpha} \wedge m^{-1}\frac{ \phi ( 4m^{1/\alpha}
  |w-y_0|) }{|w-y_0|^{d+\alpha}}  \right) .
\end{equation}
Since $m\leq
1$, the upper bound
estimate in Theorem \ref{t:main}(i)
 and \eqref{e:newqw0}--\eqref{e:newqw2} imply that
\begin{eqnarray*}
&& f_1(t, x, y)\\
  &\ge & c_2\, \int_{B(0, L_3R)^c \times B(0, L_3R)^c }
  (1\wedge \delta_D(z))^{\alpha/2}
 \left( m^{d/\alpha} \wedge  m^{-1}\frac{ \phi ( 4m^{1/\alpha}
  |x-z|) }{|x-z|^{d+\alpha}} \right)  m^{d/\alpha-d/2}
  \\
&&  \hskip 0.2truein \cdot
  (t-2/m)^{-d /2}\exp\left( -C_3 \big(m^{1/\alpha}
|z-w|\wedge m^{2/\alpha-1}\frac{|z-w|^2}{t-2/m}   \big)\right)
 (1\wedge \delta_D(w))^{\alpha/2} \\
&&  \hskip 0.2truein \cdot \left( m^{d/\alpha} \wedge m^{-1}\frac{ \phi ( 4m^{1/\alpha}
  |w-y|) }{|w-y|^{d+\alpha}} \right) dzdw \nn\\
  &\ge & c_4\,   \int_{B(0, L_3R)^c \times B(0, L_3R)^c }
 \left( m^{d/\alpha} \wedge  m^{-1}\frac{ \phi ( 4m^{1/\alpha}
  |x_0-z|) }{|x_0-z|^{d+\alpha}} \right)  m^{d/\alpha-d/2}
  \\
&&  \hskip 0.2truein \cdot (t-2/m)^{-d /2}
  \exp\left( -C_3 \big(m^{1/\alpha}
|z-w|\wedge m^{2/\alpha-1}\frac{|z-w|^2}{t-2/m}   \big)\right)
\\\
&&  \hskip 0.2truein \cdot \left( m^{d/\alpha} \wedge m^{-1}\frac{ \phi (4 m^{1/\alpha}
  |w-y_0|) }{|w-y_0|^{d+\alpha}} \right) dzdw.
   \end{eqnarray*}
Recall $A=4 \vee (C_1C_3)$.
 By the change of variables $\wh z=Az$, $\wh w=Aw$,
and  Theorem \ref{T1u}, we have that
   \begin{eqnarray}
  && f_1 (t, x, y) \nn \\
  &\ge & c_5\,    \int_{B(0, AL_3R)^c \times B(0, AL_3R)^c }
   \left( m^{d/\alpha} \wedge m^{-1} \frac{ \phi (4 A^{-1}m^{1/\alpha}
  |Ax_0-\wh z|) }{|A x_0-\wh z|^{d+\alpha}} \right) m^{d/\alpha-d/2}
   \nn\\
&&  \hskip 0.2truein \cdot
   (t-2/m)^{-d /2}
  \exp\left( -C_1^{-1}\big(m^{1/\alpha}
|\wh z-\wh w|\wedge m^{2/\alpha-1}\frac{|\wh z-\wh w|^2}{t-2/m}   \big)\right)
 \nn\\
&&  \hskip 0.2truein \cdot\left( m^{d/\alpha} \wedge  m^{-1}\frac{ \phi (4A^{-1} m^{1/\alpha}
  |\wh w-Ay_0|) }{|\wh w-Ay_0|^{d+\alpha}} \right)
  d\wh z d\wh w\nn\\
  &\ge & c_5\,    \int_{B(0, AL_3R)^c \times B(0, AL_3R)^c }
   \left( m^{d/\alpha} \wedge m^{-1} \frac{ \phi ( m^{1/\alpha}
  |Ax_0-\wh z|) }{|A x_0-\wh z|^{d+\alpha}} \right) m^{d/\alpha-d/2}
   \nn\\
&&  \hskip 0.2truein \cdot
   (t-2/m)^{-d /2}
  \exp\left( -C_1^{-1}\big(m^{1/\alpha}
|\wh z-\wh w|\wedge m^{2/\alpha-1}\frac{|\wh z-\wh w|^2}{t-2/m}   \big)\right)
 \nn\\
&&  \hskip 0.2truein \cdot\left( m^{d/\alpha} \wedge  m^{-1}\frac{ \phi ( m^{1/\alpha}
  |\wh w-Ay_0|) }{|\wh w-Ay_0|^{d+\alpha}} \right)
  d\wh z d\wh w\nn\\
      &\ge & c_6\, \int_{B(0, AL_3R)^c \times B(0, AL_3R)^c }
       p^m(1/m, A x_0,    \wh z)
  p^m(t-2/m, \wh z, \wh w)
p^m(1/m, \wh w, A y_0)
 dzdw\nn\\
  &\ge & c_6\, \int_{B(0, AL_3R)^c \times B(0, AL_3R)^c }
       p^m_{B(0, AL_3R)^c}(1/m, A x_0,    \wh z)
  p^m_{B(0, AL_3R)^c}(t-2/m, \wh z, \wh w) \nn\\
&&  \hskip 0.2truein \cdot
p^m_{B(0, AL_3R)^c}(1/m, \wh w, A y_0)
 dzdw\nn\\
      &= & c_6\,   p^m_{B(0, AL_3R)^c}(t,A x_0, A y_0).\nn
 \end{eqnarray}
 Now using \eqref{scale_kp} and  Proposition \ref{P:exterior-compare} again
 (recall that $|x_0|, |y_0| \ge L_3(L_3R)$), we conclude that
    \begin{eqnarray*}
&& p^m_{B(0, AL_3R)^c}(t,A x_0, A y_0)\,=\,
   A^{-d} p^{mA^{\alpha}}_{B(0, L_3R)^c}(A^{-\alpha} t, x_0, y_0) \\
      &\ge& c_7 m^{d/\alpha-d/2} t^{-d /2}\exp\left( -c_8 \big(A m^{1/\alpha}
|x-y|\wedge A^{2} m^{2/\alpha-1}\frac{|x_0-y_0|^2}t   \big)\right)\\
         &\ge& c_7 m^{d/\alpha-d/2} t^{-d /2}\exp\left( -c_9 \big(m^{1/\alpha}
|x-y|\wedge m^{2/\alpha-1}\frac{|x-y|^2}t   \big)\right).
  \end{eqnarray*}
Combining the last two displays with \eqref{e:lowerf1} completes the proof.
\qed

\bigskip

\noindent {\bf Proof of the lower bound estimate in Theorem \ref{t:MAIN}.}
The lower bound estimate in Theorem \ref{t:MAIN} now follows from
Theorems \ref{t:main}(i), \ref{T:exterior_1} and \ref{T:exterior}.
This completes the proof of Theorem \ref{t:MAIN}.
\qed

\section{Green function estimate}\label{S:6}

In this section, we present a  proof of Theorem \ref{T:1.3}.

\medskip

\noindent{\bf Proof of Theorem \ref{T:1.3}.}
In view of the scaling property \eqref{scale_kg} of $G^m_D$,
we may and do assume that $M=1/2$. Throughout this proof,
$m\in (0, 1/2]$.
It follows from Theorem \ref{t:MAIN} that there exists $c_i>1$, $i=1, 2,$ such that for every  $m\in (0, 1/2]$,
$t>0$ and $(x, y)\in D \times D$,
$$
p^m_D(t,x,y)\le c_1\left( 1\wedge \frac{\delta_D(x)}{ 1 \wedge t^{1/\alpha}} \right)^{\alpha/2}
\left( 1\wedge \frac{\delta_D(y)}{1 \wedge  t^{1/\alpha}}\right)^{\alpha/2}
\Psi_{d, \alpha, m, 1, c_2} ( t, x, y)
$$
and
$$
p^m_D(t,x,y)\ge c_1^{-1} \left( 1\wedge \frac{\delta_D(x)}{ 1 \wedge t^{1/\alpha}} \right)^{\alpha/2}
\left( 1\wedge \frac{\delta_D(y)}{1 \wedge  t^{1/\alpha}}
\right)^{\alpha/2} \Psi_{d, \alpha, m, 1,  1/c_2} ( t, x, y).
$$
For any $c>0$, define
$$
J_c:=\int_0^\infty \left( 1\wedge \frac{\delta_D(x)}{ 1 \wedge t^{1/\alpha}} \right)^{\alpha/2}
\left( 1\wedge \frac{\delta_D(y)}{1 \wedge  t^{1/\alpha}}\right)^{\alpha/2}
\Psi_{d, \alpha, m, 1, c} ( t, x, y) dt.
$$
Then it suffices to show that
$$
J_c\asymp \frac{1 + (m^{1/\alpha} |x-y|)^{2-\alpha}     }{|x-y|^{d-\alpha}}
 \left(1\wedge \frac{  \delta_D(x)}{ |x-y| \wedge 1}\right)^{\alpha/2}  \left(1\wedge \frac{  \delta_D(y)}{ |x-y| \wedge 1 }\right)^{\alpha/2}.
$$
Without loss of generality, we will assume $c=1$ and denote $J_1$ simply by $J$.

Using a change of variable, we see that
$$
J=I_1+\left( 1\wedge \delta_D(x)\right)^{\alpha/2}
\left( 1\wedge \delta_D(y)\right)^{\alpha/2} \left(I_2(|x-y|)+m^{d/\alpha-1}
 I_3 (m^{1/\alpha} |x-y| )\right)
$$
where
$$
I_1:=\int_{0}^1   \left( 1\wedge \frac{\delta_D(x)}{  t^{1/\alpha}} \right)^{\alpha/2}
\left( 1\wedge \frac{\delta_D(y)}{  t^{1/\alpha}}\right)^{\alpha/2}  \left( {t^{-d/\alpha}}\wedge\frac{t
  \phi(m^{1/\alpha} |x-y|)}
{|x-y|^{d+\alpha}}\right)      dt,
$$
$$
I_2(r):=\int_{1}^{1/m}    {t^{-d/\alpha}}\wedge\frac{t
  \phi(m^{1/\alpha} r)} {r^{d+\alpha}}     dt
\quad \text{and}\quad
I_3(r):= \int_{1}^\infty   t^{-d /2}
e^{- r\wedge ({r^2}/{t}) }     dt.
$$
 Note that for every $a \in [0, \infty)$
$$
 \int_{a}^{\infty}
 t^{-d/2} e^{-r^2/t} dt
=r^{2-d}  \int_0^{{r^2}/a} u^{d/2-2} e^{-u}du.
$$
Thus for $r\in (0, 1
] $
\begin{eqnarray*}
I_3(r) = r^{2-d} \int_0^{{r^2}}
u^{ d/2-2} e^{-u}du \asymp r^{2-d}\int_0^{r^2}
u^{d/2-2} du  =\frac{2}{d-2},
\end{eqnarray*}
while for $r>1$,
\begin{eqnarray*}
&&r^{2-d} \int_0^{1}
u^{d/2-2}e^{-u}du   \le \int_{r}^{\infty} t^{-d /2}    e^{ - {r^2}/{t} }     dt  \le I_3(r) \nn\\
  &&=\int_{1}^{r}
t^{-d /2}e^{-r} dt +
\int_{r}^{\infty} t^{-d /2} e^{ - {r^2}/{t} }      dt
 \,\le\, c_3 \, e^{-r}     +r^{2-d} \int_0^{r}
u^{d/2-2}e^{-u}du  \,\le\,  c_4 \, r^{2-d}.
\end{eqnarray*}
Thus we have
\begin{equation}\label{e:G1}
I_3 (r) \asymp 1\wedge  r^{2-d}.
\end{equation}
Noting that $m\in (0, 1/2]$, so when $m^{1/\alpha} r \le 1$,
\begin{eqnarray} \label{e:G2}
 I_2 (r) &\asymp&
 \int_{1}
 ^{1/m}
  \left( {t^{-d/\alpha}}\wedge\frac{t}
  {r^{d+\alpha}}\right) dt
  = \int_1^{r^\alpha \vee 1} \frac{t}{r^{d+\alpha}}dt  +
    \int_{r^\alpha \vee 1}^{1/m} t^{-d/\alpha} dt \nn \\
  &=& \frac{1}{2r^{d+\alpha}} \left( (r^\alpha \vee 1)^2-1\right)
  + \frac{\alpha}{d-\alpha} \left( (r^\alpha \vee 1)^{1-d/\alpha}
  - m^{d/\alpha -1} \right) \nn \\
  &\asymp & 1\wedge r^{\alpha -d}.
\end{eqnarray}
If $m^{1/\alpha} r >1$, then
$$
  I_2 (r)\asymp \int_{1}^{1/m}     \left( {t^{-d/\alpha}}\wedge\frac{t
e^{-m^{1/\alpha} r}  }
{r^{d+\alpha}}\right)      dt
$$
and  a change of variable $s=tm$ gives
\begin{eqnarray}
I_2(r) &\asymp&  m^{d/\alpha -1} \int_m^1
\left( {s^{-d/\alpha}}\wedge\frac{s e^{-m^{1/\alpha} r}}
{(m^{1/\alpha} r)^{d+\alpha}}\right)
ds = m^{d/\alpha-1} \int_m^1  \frac{s e^{-m^{1/\alpha} r} }
{(m^{1/\alpha} r)^{d+\alpha}} ds \nn \\
&\asymp&   m^{d/\alpha -1}\frac{ e^{-m^{1/\alpha} r}  }
{(m^{1/\alpha} r)^{d+\alpha}}
=\frac{e^{-m^{1/\alpha} r}}{ (m^{1/\alpha} r)^2r^{d-\alpha}}.
\label{e:G3}
\end{eqnarray}

\medskip \noindent
 (i) Suppose $m^{1/\alpha} |x-y| \le 1$. Since $\phi(m^{1/\alpha} |x-y|) \asymp 1$, it follows from \cite[(4.3), (4.4) and (4.6)]{CKS} that
$$
I_1 \asymp \frac{1} {|x-y|^{d-\alpha}}
 \left(1\wedge \frac{  \delta_D(x)}{ |x-y| }\right)^{\alpha/2}  \left(1\wedge \frac{  \delta_D(y)}{ |x-y| }\right)^{\alpha/2}.
$$
Thus, we have by \eqref{e:G1} and  \eqref{e:G2} that
\begin{eqnarray}\label{e:G4}
J &\asymp &  \frac{1}{|x-y|^{d-\alpha}}
  \left(1\wedge \frac{  \delta_D(x)}{ |x-y| }\right)^{\alpha/2}  \left(1\wedge \frac{  \delta_D(y)}{ |x-y| }\right)^{\alpha/2} \nn \\
  && + \left( 1\wedge \delta_D(x)\right)^{\alpha/2}
\left( 1\wedge \delta_D(y)\right)^{\alpha/2} \left( 1\wedge \frac{1}{|x-y|^{d-\alpha}} + m^{d/\alpha -1} \right)\nn\\
&\asymp &   \frac{1} {|x-y|^{d-\alpha}}
 \left(1\wedge \frac{  \delta_D(x)}{ |x-y| \wedge 1}\right)^{\alpha/2}  \left(1\wedge \frac{  \delta_D(y)}{ |x-y| \wedge 1 }\right)^{\alpha/2}.
\end{eqnarray}
We have arrived the last display above by considering the cases $|x-y|\geq 1$
and $|x-y|< 1$ separately.

\medskip \noindent
 (ii) Suppose $m^{1/\alpha} |x-y| > 1$.
For $0 \le t \le 1$, we have $0 \le mt \le 1/2 $ and so
$$
\left( t^{-d/\alpha}
\wedge \frac{t\phi(m^{1/\alpha} |x-y|)}{|x-y|^{d+\alpha}}\right)=
m^{d/\alpha} \left( (mt)^{-d/\alpha}
\wedge \frac{(mt)\phi(m^{1/\alpha} |x-y|)}{(m^{1/\alpha}|x-y|)^{d+\alpha}}\right) =\frac{t\phi(m^{1/\alpha} |x-y|)}{|x-y|^{d+\alpha}} .
$$
Thus  by the change of variable $u= \frac{|x-y|^\alpha}{t}$, we have
\begin{eqnarray}
I_1&=& \frac{\phi (m^{1/\alpha} |x-y|)}{|x-y|^{d+\alpha}}
\int_0^1 t \left( 1\wedge \frac{\delta_D(x)}{t^{1/\alpha}}\right)^{\alpha/2}
 \left( 1\wedge \frac{\delta_D(y)}{t^{1/\alpha}}\right)^{\alpha/2} dt \nn \\
&=&\frac{\phi(m^{1/\alpha} |x-y|) }{|x-y|^{d-\alpha}}  \int_{|x-y|^\alpha}^\infty
 u^{-3} \left(1\wedge
\frac{ {\sqrt u} \delta_D(x)^{\alpha/2}}{ |x-y|^{\alpha/2} }\right)
\left(1\wedge \frac{ {\sqrt u} \delta_D(y)^{\alpha/2}}{
|x-y|^{\alpha/2} }\right) du.          \label{ID:5.1}
\end{eqnarray}
Note that  since $|x-y| \ge m^{-1/\alpha} > 2^{1/\alpha}$,
\begin{eqnarray}
&&   \int_{|x-y|^\alpha}^\infty
 u^{-3}  \left(1\wedge
\frac{ {\sqrt u} \delta_D(x)^{\alpha/2}}{ |x-y|^{\alpha/2} }\right)
\left(1\wedge \frac{ {\sqrt u} \delta_D(y)^{\alpha/2}}{
|x-y|^{\alpha/2} }\right) du   \nonumber\\
&=&
   \int_{|x-y|^\alpha}^\infty
  u^{-2} \, \left(u^{-1/2} \wedge
\frac{   \delta_D(x)^{\alpha/2}}{ |x-y|^{\alpha/2} }\right)
\left(u^{-1/2} \wedge \frac{  \delta_D(y)^{\alpha/2}}{
|x-y|^{\alpha/2} }\right) du \nonumber\\
&\leq&
 \int_{|x-y|^\alpha}^\infty
  u^{-2} \,
  du\left(1 \wedge
\frac{   \delta_D(x)^{\alpha/2}}{ |x-y|^{\alpha/2} }\right) \left(1
\wedge \frac{  \delta_D(y)^{\alpha/2}}{
|x-y|^{\alpha/2} }\right)  \nonumber\\
&=& |x-y|^{-\alpha}    \left(1\wedge \frac{
\delta_D(x)}{ |x-y| }\right)^{\alpha/2} \left(1\wedge
\frac{  \delta_D(y)}{ |x-y|  }\right)^{\alpha/2}.
\nonumber
\end{eqnarray}
Thus it follows from  \eqref{ID:5.1} that
\begin{eqnarray}\label{e:G7}
I_1 &\leq&  \frac{\phi (1)}{|x-y|^{d}} \left(1\wedge \frac{
\delta_D(x)}{ |x-y| }\right)^{\alpha/2} \left(1\wedge
\frac{  \delta_D(y)}{ |x-y|  }\right)^{\alpha/2} \nn \\
&\leq &   \frac{\phi (1) \, m^{2/\alpha}} {|x-y|^{d-2}}
\left(1\wedge \delta_D(x) \right)^{\alpha/2}
\left(1\wedge  \delta_D(y) \right)^{\alpha/2},
\end{eqnarray}
where in the last inequality, we used the assumption
 $m^{1/\alpha}|x-y|\geq 1$. Recalling $m\in (0, 1/2]$,
we thus have by \eqref{e:G1}, \eqref{e:G3} and \eqref{e:G7} that
\begin{eqnarray}
J&\asymp& \left(1\wedge \delta_D(x) \right)^{\alpha/2}
\left(1\wedge  \delta_D(y) \right)^{\alpha/2}
\left( \frac{e^{-m^{1/\alpha} |x-y|}}{ (m^{1/\alpha} |x-y|)^2|x-y|^{d-\alpha}} +   \frac{m^{2/\alpha-1}}{|x-y|^{d-2}} \right) \nn\\
&=& \left(1\wedge \delta_D(x) \right)^{\alpha/2}
\left(1\wedge  \delta_D(y) \right)^{\alpha/2}
\left( \frac{e^{-m^{1/\alpha} |x-y|}}{ (m^{1/\alpha} |x-y|)^2|x-y|^{d-\alpha}} +   \frac{(m^{1/\alpha} |x-y|)
^{2-\alpha}}{|x-y|^{d-\alpha}} \right) \nn\\
&\asymp & \left(1\wedge \delta_D(x) \right)^{\alpha/2}
\left(1\wedge  \delta_D(y) \right)^{\alpha/2}
   \frac{(m^{1/\alpha} |x-y|)
^{2-\alpha}}{|x-y|^{d-\alpha}}  \nn \\
&\asymp &   \frac{1+(m^{1/\alpha} |x-y|)
^{2-\alpha}}{|x-y|^{d-\alpha}}  \left(1\wedge \frac{  \delta_D(x)}{ |x-y| \wedge 1}\right)^{\alpha/2}  \left(1\wedge \frac{  \delta_D(y)}{ |x-y| \wedge 1 }\right)^{\alpha/2}. \nn
\end{eqnarray}
This combing with \eqref{e:G4} completes the proof of the theorem
\qed

\medskip

{\bf Zhen-Qing Chen}

Department of Mathematics, University of Washington, Seattle,
WA 98195, USA

E-mail: \texttt{zqchen@uw.edu}

\bigskip

{\bf Panki Kim}

Department of Mathematics, Seoul National University,
Building 27, 1 Gwanak-ro, Gwanak-gu
Seoul 151-747, Republic of Korea

E-mail: \texttt{pkim@snu.ac.kr}

\bigskip

{\bf Renming Song}

Department of Mathematics, University of Illinois, Urbana, IL 61801, USA

E-mail: \texttt{rsong@math.uiuc.edu}
\end{document}